\documentclass[reqno]{amsart}
\usepackage{graphicx,amsfonts,amssymb,amsmath,amsthm,epsfig,euscript}

\newtheorem{theorem}{Theorem}
\newtheorem{corollary}[theorem]{Corollary}
\newtheorem{proposition}[theorem]{Proposition}

\newtheorem{definition}[theorem]{Definition}
\newtheorem{example}[theorem]{Example}
\newtheorem{remark}[theorem]{Remark}

\newcommand{\sg}{\sigma}
\newcommand{\fig}[3]
{\begin{figure}[ht] \centerline{\scalebox{#1}{\epsfig{file=#2.eps}}}
\vspace{-1mm} \caption{#3 \label{fig:#2}}
\end{figure}}

\numberwithin{theorem}{section}
\numberwithin{equation}{section}

\begin{document}
\title{Counting descent pairs with prescribed tops and bottoms}

\author{John T. Hall}
\author{Jeffrey B. Remmel}\thanks{Partially supported by NSF grant DMS 0400507}
\address{Department of Mathematics\\
University of California, San Diego\\
La Jolla, CA 92093}
\email{jthall@math.ucsd.edu}
\email{jremmel@ucsd.edu}

\keywords{permutation patterns, descents, excedences, descent tops, descent bottoms, rook placements}
\subjclass[2000]{05A05, 05A15}

\begin{abstract}
Given sets $X$ and $Y$ of positive integers and a permutation $\sigma = \sigma_1 \sigma_2 \cdots \sigma_n \in 
S_n$, an 
$X,Y$-\emph{descent} of $\sigma$ is a descent pair $\sigma_i > \sigma_{i+1}$ whose ``top'' $\sigma_i$ is in 
$X$ and whose ``bottom'' $\sigma_{i+1}$ is in $Y$. 
We give two formulas for the number $P_{n,s}^{X,Y}$ of $\sigma \in S_n$ with $s$ $X,Y$-descents.
$P_{n,s}^{X,Y}$ is also shown to be a hit number of a certain Ferrers board. This work generalizes results of Kitaev and 
Remmel \cite{KitaevRemmel1} \cite{KitaevRemmel2} on counting descent pairs whose top (or bottom) is 
equal to $0 \mod k$.
\end{abstract}

\maketitle

\tableofcontents

\section{Introduction}
\label{introduction-section}

Let $S_n$ denote the set of permutations of the set $[n]= \{1, 2, \ldots, n\}$. 
A \emph{descent pair} of a permutation $\sigma = \sigma_1 \sigma_2 \cdots \sigma_n \in S_n$ is a 
pair $(\sigma_i, \sigma_{i+1})$ with $\sigma_i > \sigma_{i+1}$.  The main focus of this paper is to 
study the distribution of descent pairs whose top $\sigma_i$ lies in some fixed 
set $X$ and whose bottom $\sigma_{i+1}$ lies in some fixed set $Y$.

\begin{definition}
Given subsets $X,Y \subseteq \mathbb{N}$ and a permutation $\sigma \in S_n$, let
\begin{eqnarray*}
Des_{X,Y}(\sigma) & = &  \{ i : 
\sigma_i > \sigma_{i+1} \ \& \ \sigma_i \in X \ \& \ \sigma_{i+1} \in Y \}, \mbox{ and} \\
des_{X,Y}(\sigma) & = & 
|Des_{X,Y}(\sigma)|.
\end{eqnarray*}
If $i \in Des_{X,Y}(\sg)$, then we call the pair $(\sg_i,\sg_{i+1})$ an 
$X,Y$-\emph{descent}. 
\end{definition}
For example, if $X=\{2,3,5\}, Y = \{1,3, 4\},$ and $\sigma = 54213$, then $\mbox{\emph{Des}}_{X,Y}(\sigma) = 
\{ 1,3 \}$ and $\mbox{\emph{des}}_{X,Y}(\sigma) = 2$.

For fixed $n$ we define the polynomial
\begin{equation}\label{Pdef}
P_n^{X,Y}(x)  = \sum_{s \geq 0} P_{n,s}^{X,Y}x^s := \sum_{\sg \in S_n} x^{des_{X,Y}(\sg)}.
\end{equation}
Thus the coefficient $P_{n,s}^{X,Y}$ is the number of $\sigma \in S_n$ with exactly $s$ $X,Y$-descents. 

Our main result is to give direct combinatorial proofs of a pair of formulas for $P_{n,s}^{X,Y}$. First of all, 
for any set $S \subseteq \mathbb{N}$, let
\begin{eqnarray*}
S_n & = & S \cap [n], \mbox{ and} \\
S_n^c & = & (S^c)_n = [n]-S.
\end{eqnarray*}
Then we have

\noindent {\bf Theorem 2.3}
\begin{equation}\label{Ieq:combXY1}
P_{n,s}^{X,Y} = \left| X_n^c\right|!\sum_{r=0}^s (-1)^{s-r} \binom{ \left| X_n^c\right| +r}{r}
\binom{n+1}{s-r}
\prod\limits_{x \in X_n} (1+r + \alpha_{X,n,x} + \beta_{Y,n,x}),
\end{equation}
and

\noindent {\bf Theorem 2.5}
\begin{equation}\label{Ieq:combXY2}
P_{n,s}^{X,Y} = \left| X_n^c\right| !\sum_{r=0}^{ \left| X_n\right|-s} (-1)^{\left| X_n\right|-s-r} 
\binom{\left| X_n^c\right|+r}{r}
\binom{n+1}{ \left| X_n \right|-s-r}
\prod\limits_{x \in X_n} (r+\beta_{X,n,x}-\beta_{Y,n,x}),
\end{equation}
where for any set $S$ and any $j,1 \leq j \leq n$, we define  
\begin{eqnarray*}
\alpha_{S,n,j} &=&  |S^c \cap \{j+1, j+2,\ldots, n\}| = |\{x: j < x \leq n \ \& \ x \notin S\}|,\mbox{ and} \\
\beta_{S,n,j} &=&  |S^c \cap \{1, 2, \ldots, j-1\}| = |\{x: 1\leq  x < j \ \& \ x \notin S\}|.
\end{eqnarray*}

\begin{example}
Suppose $X=\{2,3,4,6,7,9\}, Y = \{1, 4 , 8\}$, and $n=6$. Thus $X_6 = \{2, 3, 4, 6 \},X_6^c = \{1, 5 
\}, Y_6 = \{1, 4 \}, Y_6^c = \{ 2, 3, 5, 6 \}$, and we have the following table of values of 
$\alpha_{X,6,x}, \beta_{Y,6,x}$, and $\beta_{X,6,x}$.
\begin{center}
\begin{tabular}{|c|c|c|c|c|}
\hline
$x$ & $2$ & $3$ & $4$ & $6$ \\
\hline
$\alpha_{X,6,x}$ & $1$ & $1$ & $1$ & $0$ \\
\hline
$\beta_{Y,6,x}$ & $0$ & $1$ & $2$ & $2$ \\
\hline
$\beta_{X,6,x}$ & $1$ & $1$ & $1$ & $2$ \\
\hline
\end{tabular}
\end{center}
(\ref{Ieq:combXY1}) gives
\begin{eqnarray*}
P_{6,2}^{X,Y} & = & 2!\sum\limits_{r=0}^2 (-1)^{2-r} \binom{2+r}{r} \binom{7}{2-r} (2+r)(3+r)(4+r)(3+r) \\
 & = & 2\left(1 \cdot 21 \cdot 2 \cdot 3 \cdot 4 \cdot 3 - 3 \cdot 7 \cdot 3 \cdot 4 \cdot 5 \cdot 4 + 6 
\cdot 1 \cdot 4 \cdot 5 \cdot 6 \cdot 5 \right) 
\\
 & = & 2(1512 - 5040 + 3600) \\
 & = & 144.
\end{eqnarray*}
while (\ref{Ieq:combXY2}) gives
\begin{eqnarray*}
P_{6,2}^{X,Y} & = & 2!\sum\limits_{r=0}^2 (-1)^{2-r} \binom{2+r}{r} \binom{7}{2-r} (1+r)(0+r)(-1+r)(0+r) \\
 & = & 2 \left( 1 \cdot 21 \cdot 1 \cdot 0 \cdot (-1) \cdot 0 - 3 \cdot 7 \cdot 2 \cdot 1 \cdot 0 \cdot 1 
 + 6 \cdot 1 \cdot 3 \cdot 2 \cdot 1 \cdot 2 \right) \\
 & = & 2(0 - 0 + 72) \\
 & = & 144.
\end{eqnarray*}
\end{example}

Let $P_{n,s}^X := P_{n,s}^{X,\mathbb{N}}$. Since $\beta_{\mathbb{N},n,x} = 0$ for all $x \in X_n$, we have as 
corollaries
\begin{equation}\label{eq:fund1}
P_{n,s}^X = \left| X_n^c\right| !\sum_{r=0}^s (-1)^{s-r} \binom{\left| X_n^c\right|+r}{r} \binom{n+1}{s-r} 
\prod_{x \in X_n} (1+r + \alpha_{X,n,x})
\end{equation} 
and 
\begin{equation}\label{eq:fund2}
P_{n,s}^X = \left| X_n^c\right| !\sum_{r=0}^{\left| X_n\right|-s} (-1)^{\left| X_n\right|-s-r}
\binom{\left| X_n^c\right|+r}{r} \binom{n+1}{\left| X_n \right|-s-r} 
\prod_{x \in X_n} (r + \beta_{X,n,x}).
\end{equation} 

We will show that the equality of (\ref{eq:fund1}) and (\ref{eq:fund2}) is 
equivalent to a certain special case of a general transformation result due 
to Gasper \cite{Gasper} for hypergeometric series of Karlsson-Minton type. Since (\ref{eq:fund1}) and (\ref{eq:fund2}) have completely combinatorial proofs, it follows that we can give combinatorial proofs of many special cases of Gasper's transformation theorem. 

We will use the fundamental transformation of Foata \cite{FoataSchutzenberger1} to show 
that the polynomials $P_n^{X,Y}(x)$ are special cases of hit polynomials for Ferrers boards $B$ contained in the $n 
\times 
n$ board. Since there are many formulas for hit polynomials of Ferrers boards 
(see, for example, the results of Haglund \cite{Haglund}), we have an alternative way to prove formulas 
for the polynomials $P_n^{X,Y}(x)$. Moreover, we will show that we can use the same idea to reduce the computation of 
the coefficients of  $P_n^{X,Y}(x)$ to the computation of the coefficients of  
$P_n^{X^*,\mathbb{N}}(x)$, for some appropriate $X^*$ depending on $X$ and $Y$. Thus (\ref{eq:fund1}) and 
(\ref{eq:fund2}) already contain all of the information needed to compute the seemingly more general formulas 
(\ref{Ieq:combXY1}) and (\ref{Ieq:combXY2}).

This type of study was initiated by 
Kitaev and Remmel \cite{KitaevRemmel1,KitaevRemmel2}.
In particular, they studied descents according to the equivalence class mod $k$ of either
the top or bottom of a descent pair. For any
set $X \subseteq \{1,2,3, \ldots \}$, they defined
\begin{itemize}
\item $\overleftarrow{Des}_X(\sg) = \{i: \sg_i > \sg_{i+1} \ \& \
\sg_i \in X\}$ and $\overleftarrow{des}_X(\sg) =
|\overleftarrow{Des}_X(\sg)|$, and 
\item $\overrightarrow{Des}_X(\sg) = \{i: \sg_i > \sg_{i+1} \ \& \
\sg_{i+1} \in X\}$ and $\overrightarrow{des}_X(\sg) =
|\overrightarrow{Des}_X(\sg)|$.
\end{itemize}
It is easy to see that $\overleftarrow{Des}_X(\sigma) = 
Des_{X,\mathbb{N}}(\sigma)$ and $\overrightarrow{Des}_X(\sigma) = 
Des_{\mathbb{N},X}(\sigma)$.
In \cite{KitaevRemmel1}, Kitaev and Remmel studied polynomials such as
\begin{eqnarray*}
R_n(x) & = & \sum_{k \geq 0} R_{k,n} x^k := \sum_{\sg \in S_n} x^{\overleftarrow{des}_E(\sg)}, \mbox{ and}\\
Q_n(x) & = & \sum_{k \geq 0} Q_{k,n} x^k := \sum_{\sg \in S_n} x^{\overrightarrow{des}_E(\sg)}, \\
\end{eqnarray*}
where $E = \{2,4,6,  \ldots, \}$ is the set of positive even integers. In these cases, 
they found surprisingly simple formulas for the coefficients. For example, they showed that
\begin{equation}\label{eq:KR1}
R_{2n,k} = \binom{n}{k}^2 (n!)^2.
\end{equation}
In \cite{KitaevRemmel2}, Kitaev and Remmel studied polynomials such as
\begin{eqnarray*}
A^{(k)}_n(x) & = & \sum_{j \geq 0} A^{(k)}_{n,j} x^j := \sum_{\sg \in S_n}
x^{\overleftarrow{des}_{kN}(\sg)}, \mbox{ and}\\
B^{(k)}_n(x) & = & \sum_{j \geq 0} B^{(k)}_{n,j} x^j := \sum_{\sg \in S_n}
x^{\overrightarrow{des}_{kN}(\sg)},
\end{eqnarray*}
where $kN = \{k,2k,3k, \ldots \}$. Note that $A_n^{(k)}(x) = P_n^{X,\mathbb{N}}(x)$ and 
$B_n^{(k)}(x) =P_n^{\mathbb{N},Y}(x),$ where $X$ and $Y$ are $k\mathbb{N}$ for some $k \geq 2$. 
When $k \geq 3$, the formulas for 
the coefficients of these polynomials are not as simple as (\ref{eq:KR1}). 
For example, Kitaev and Remmel used a recursion to prove the following formulas, which hold for all $0 \leq 
j \leq k-1$ and all $n \geq 0$.
\begin{eqnarray}\label{eq:KR2}
A^{(k)}_{kn+j,s} &=& ((k-1)n+j)!  \times \nonumber \\
&& \ \ \sum_{r=0}^s (-1)^{s-r}\binom{(k-1)n+j+r}{r}
\binom{kn+j+1}{s-r} \prod_{i=0}^{n-1} (r+1+j +(k-1)i)
\end{eqnarray}
\begin{eqnarray}\label{eq:KR3}
A^{(k)}_{kn+j,s} &=& ((k-1)n+j)!  \times \nonumber \\
&& \ \ \sum_{r=0}^{n-s} (-1)^{n-s-r}\binom{(k-1)n+j+r}{r}
\binom{kn+j+1}{n-s-r} \prod_{i=1}^{n} (r+(k-1)i)
\end{eqnarray}
Our main results are generalizations of the two formulas (\ref{eq:KR2}) and (\ref{eq:KR3}). 

The outline of this paper is as follows. 
In Section 2, we give several
formulations of a recursion for the number $P_{n,s}^{X,Y}$ 
of $\sigma \in S_n$ with $s$ $X,Y$-descents. Then we present our main results, the combinatorial
proofs of (\ref{Ieq:combXY1}) and (\ref{Ieq:combXY2}). 
In Section 3,  we give several applications of our main results, including new proofs of results of Kitaev and 
Remmel \cite{KitaevRemmel1} \cite{KitaevRemmel2} on counting descent pairs whose top (or bottom) is
equivalent to $0 \mod k$, and  a combinatorial proof of various special cases of a transformation of 
Karlsson-Minton type hypergeometric series due to Gasper \cite{Gasper}. In Section 4, we use Foata's First 
Transformation, which is a bijection taking descents to excedences, to rephrase the problem of computing the 
polynomials $P^{X,Y}_n(x)$ as one of computing hit polynomials for certain boards 
$B$ contained in the $n \times n$ board. 
In Section 5, we show that our results can be extended from permutations to words. In Section 6, we 
shall consider a more general problem. That is, 
for any $X,Y,Z \subseteq \mathbb{N}$, we can consider the polynomials 
\begin{equation*}
P_n^{X,Y,Z}(x) = \sum_{s \geq 0} P_{n,s}^{X,Y,Z}x^s := \sum_{\sg \in S_n} x^{des_{X,Y,Z}(\sg)},
\end{equation*}
where for any subsets $X$, $Y$, and $Z$ of $\mathbb{N}$, and 
permutation $\sg = \sg_1 \sg_2 \cdots \sg_n \in S_n$,
\begin{eqnarray*}
Des_{X,Y,Z}(\sg) & = & \{i:\sg_i > \sg_{i+1} \ \& \ \sg_i \in X, \sg_{i+1} \in 
Y \ \& \ \sg_i - \sg_{i+1} \in Z\}, \mbox{ and} \\
des_{X,Y,Z}(\sg) & = & |Des_{X,Y,Z}(\sg)|.
\end{eqnarray*}
Clearly $P_n^{X,Y}(x) = P_n^{X,Y,\mathbb{N}}(x)$. We do not have a  
formula for the coefficients $P_{n,s}^{X,Y,Z}$ for arbitrary $X$, $Y$, and 
$Z$. However, we will show that we can find formulas for 
$P_{n,s}^{X,Y,Z}$ for certain special cases of $X$, $Y$, and $Z$.
Finally, in Section 7, we present 
some directions for future research.

\section{Prescribed Tops and Bottoms}
\label{topbottom-section}
In this section, we will give several ways to compute the coefficients 
$P_{n,s}^{X,Y}$. 

Given $X,Y \subseteq \mathbb{N}$, let $P_0^{X,Y}(x,y) = 1$, and for $n \geq 1$, define 
\begin{equation*}\label{eq:PXY} 
P_n^{X,Y}(x,y) = 
\sum_{s,t \geq 0} P_{n,s,t}^{X,Y} x^sy^t :=\sum\limits_{\sigma 
\in S_n} 
x^{des_{X,Y}(\sigma)}y^{|Y_n^c|}.
\end{equation*}
Let 
$\Phi_{n+1}$ and $\Psi_{n+1}$ be the operators defined as
\begin{equation*}
\begin{array}{l}
\Phi_{n+1} : x^s y^t \longrightarrow sx^{s-1} y^t + (n+1-s) x^sy^t \\
\Psi_{n+1} : x^s y^t \longrightarrow (s+t+1) x^s y^t + (n-s-t) x^{s+1}y^t.
\end{array}
\end{equation*}

\begin{proposition}\label{insert3}
For any sets $X,Y \subseteq \mathbb{N}$, the polynomials $P_n^{X,Y}(x,y)$ satisfy
\begin{equation*}
P_{n+1}^{X,Y}(x,y)  = \left\{ \begin{array}{ll} y\cdot\Phi_{n+1} (P_n^{X,Y}(x,y)) & \mbox{ if } 
n+1 \not\in X \mbox{ and }  n+1 \not\in Y, \\
\Phi_{n+1} (P_n^{X,Y}(x,y)) & \mbox{ if } 
n+1 \not\in X \mbox{ and } n+1 \in Y, \\
y\cdot \Psi_{n+1} (P_n^{X,Y}(x,y)) & \mbox{ if } 
n+1 \in X \mbox{ and }  n+1 \not\in Y, \mbox{ and}\\
\Psi_{n+1} (P_n^{X,Y}(x,y)) & \mbox{ if } 
n+1 \in X \mbox{ and }  n+1 \in Y. \end{array} \right.
\end{equation*}
\end{proposition}

\begin{proof}
We think of a permutation as being built up by successively inserting the numbers
$1, 2, 3$, and so on. Given a permutation $\sigma \in S_n$ with $s$ $X,Y$-descents, if $n+1 \not\in X$
and we insert $n+1$ in the middle of one of the $s$ $X,Y$-descent pairs, then
we destroy that $X,Y$-descent and get a permutation with $s-1$ $X,Y$-descents. 
If, instead, we insert $n+1$ in one of the other $n+1-s$ possible spots (including the spots at the beginning and end of the permutation), 
then we preserve the number of $X,Y$-descents. Thus, if $n+1 \notin X$, we have
\begin{equation*}
P_{n+1}^{X,Y}(x,y)  = \left\{ \begin{array}{ll} y\cdot\Phi_{n+1} (P_n^{X,Y}(x,y)) & \mbox{ if } 
n+1 \not\in Y, \mbox{ and } \\
\Phi_{n+1} (P_n^{X,Y}(x,y)) & \mbox{ if } 
 n+1 \in Y. \end{array} \right.
\end{equation*}
On the other hand, if $n+1 \in X$, then we 
preserve the number of $X,Y$-descents by inserting $n+1$ in the middle of one of the $s$
$X,Y$-descent pairs, or before any of the $t$ elements of $Y_n^c$, or at the end of the permutation.
If, instead, we insert $n+1$ in one of the 
other $n-s-t$ possible spots we create a new $X,Y$-descent. Thus, if $n+1 \in X$, we have
\begin{equation*}
P_{n+1}^{X,Y}(x,y)  = \left\{ \begin{array}{ll} y\cdot\Psi_{n+1} (P_n^{X,Y}(x,y)) & \mbox{ if } 
n+1 \not\in Y, \mbox{ and } \\
\Psi_{n+1} (P_n^{X,Y}(x,y)) & \mbox{ if } 
 n+1 \in Y. \end{array} \right.
\end{equation*}
\end{proof}

It is easy to see that Proposition \ref{insert3} implies the following result.
\begin{corollary}
For all $X,Y \subseteq \mathbb{N}$ and $n \geq 1$, the following recursion 
holds for the coefficients $P^{X,Y}_{n,s,t}$.
\begin{equation}\label{XYrec}
P_{n+1,s,t}^{X,Y} = \left\{ \begin{array}{ll} (s+1)P_{n,s+1,t-1}^{X,Y}+(n+1-s)P_{n,s,t-1}^{X,Y} & 
\mbox{ if } n+1 \not\in X \mbox{ and } n+1 \not\in Y,\\
(s+1)P_{n,s+1,t}^{X,Y}+(n+1-s)P_{n,s,t}^{X,Y} &
\mbox{ if } n+1 \not\in X \mbox{ and } n+1 \in Y,\\
(s+t)P_{n,s,t-1}^{X,Y}+(n+2-s-t)P_{n,s-1,t-1}^{X,Y} &
\mbox{ if } n+1 \in X \mbox{ and } n+1 \not\in Y, \mbox{ and}\\
(s+t+1)P_{n,s,t}^{X,Y}+(n+1-s-t)P_{n,s-1,t}^{X,Y} & 
\mbox{ if } n+1 \in X \mbox{ and } n+1 \in Y. \end{array} \right.
\end{equation}
\end{corollary}

We can rephrase Proposition \ref{insert3} in terms of partial differential equations as follows.
\tiny
\begin{equation}\label{diffeq3}
P_{n+1}^{X,Y}(x,y) = \left\{ \begin{array}{ll} y\left( (n+1)P_n^{X,Y}(x,y) + (1-x)\frac{\partial}{\partial x}P_n^{X,Y}(x,y) \right) & 
\mbox{ if } n+1 \not\in X \mbox{ and } n+1 \not\in Y, \\
(n+1)P_n^{X,Y}(x,y) + (1-x)\frac{\partial}{\partial x}P_n^{X,Y}(x,y) & 
\mbox{ if } n+1 \not\in X \mbox{ and } n+1 \in Y, \\
y\left( (xn+1)P_n^{X,Y}(x,y) + (x-x^2)\frac{\partial}{\partial x}P_n^{X,Y}(x,y) + 
y(1-x)\frac{\partial}{\partial y}P_n^{X,Y}(x,y) \right) &
\mbox{ if } n+1 \in X \mbox{ and } n+1 \not\in Y, \mbox{ and} \\
(xn+1)P_n^{X,Y}(x,y) + (x-x^2)\frac{\partial}{\partial x}P_n^{X,Y}(x,y) + 
y(1-x)\frac{\partial}{\partial y}P_n^{X,Y}(x,y) & 
\mbox{ if } n+1 \in X \mbox{ and } n+1 \in Y. \end{array} \right.
\end{equation}
\normalsize

One can then use either of the recursions (\ref{XYrec}) and (\ref{diffeq3}) to 
compute $P_n^{X,Y}(x,y)$ for any $X$ and $Y$. For example, if 
$X = \{2,3,5 \}$ and $Y= \{1, 3, 4\}$ we have
\begin{eqnarray*}
P_0^{X,Y}(x,y) & = & 1 \\
P_1^{X,Y}(x,y) & = & 1 \\
P_2^{X,Y}(x,y) & = & y(1+x)\\
P_3^{X,Y}(x,y) & = & y(2+4x) \\
P_4^{X,Y}(x,y) & = & y(12 + 12x) \\
P_5^{X,Y}(x,y) & = & y^2(24+72x+24x^2).
\end{eqnarray*}

It is easy to see that $P_{n,s,t}^{X,Y} = 0$ unless $t = \left| Y_n^c \right|$, so we will drop
the $t$ and write $P_{n,s}^{X,Y}$ for $P_{n,s,|Y^c_n|}^{X,Y}$.
Thus $P_n^{X,Y}(x,1)$ is really just the polynomial $P_n^{X,Y}(x)$ defined by (\ref{Pdef}).

Recall that for each $x \in X_n$, we defined
\begin{eqnarray*}
\alpha_{X,n,x} &=& |\{z: x < z \leq n \ \& \ z \notin X\}|, \mbox{ and}\\ 
\beta_{Y,n,x} &=& |\{z: 1\leq  z < x \ \& \ z \notin Y\}|.
\end{eqnarray*}
Then we have following formula for $P_{n,s}^{X,Y}$.

\begin{theorem}\label{combXY1}
\begin{equation}\label{eq:combXY1}
P_{n,s}^{X,Y} = \left| X_n^c \right| !\sum_{r=0}^s (-1)^{s-r} \binom{ \left| X_n^c \right|+r}{r}
\binom{n+1}{s-r}
\prod\limits_{x \in X_n} (1+r +\alpha_{X,n,x} + \beta_{Y,n,x}).
\end{equation}
\end{theorem}

\begin{remark}
Theorem \ref{combXY1} can be proved by showing that the formula satisfies the recursion (\ref{XYrec}) for $P^{X,Y}_{n,s}$. 
However, we will give a direct combinatorial proof using a sign-reversing involution on a set of \emph{configurations}, which are arrays 
of numbers, $+$'s, and $-$'s. The basic idea is simple: applying the involution to each configuration results in either changing a $+$ to a $-$, 
or changing a $-$ to a $+$. The fixed points of the involution
will be shown to correspond naturally to permutations $\sigma \in S_n$ such that $des_{X,Y}(\sg) = s$.
\end{remark}

\begin{proof} Let $X,Y,n$, and $s$ be given. For $r$ satisfying $0 \leq r \leq s$, we define the 
set of what we call $(n,s,r)^{X,Y}$-\emph{configurations}. An 
$(n,s,r)^{X,Y}$-configuration $c$ consists of an array of the 
numbers $1, 2, \ldots, n$, $r$ $+$'s, and $(s-r)$ $-$'s, satisfying the following two conditions:
\begin{enumerate}
\item[(i)] each $-$ is either at the very beginning of the array or immediately follows a number, and
\item[(ii)] if $x \in X$ and $y\in Y$ are consecutive numbers in the array, and $x>y$, i.e., if $(x,y)$ 
forms an $X,Y$-descent pair in the underlying permutation, then there must be at least one $+$ between $x$ and $y$.
\end{enumerate}
Note that in an $(n,s,r)^{X,Y}$-configuration, the number of $+$'s plus the number of $-$'s equals $s$.
For example, if $X = \{2, 3, 5, 6 \}$  and $Y = \{1, 3\}$, the following is a $(6,5,3)^{X,Y}$-configuration.
\begin{equation*}
c = 5+2-+46+13-
\end{equation*}
In this example, the underlying permutation is $524613$. In general, we will let $c_1c_2 \cdots c_n$ 
denote the underlying permutation of the $(n,s,r)^{X,Y}$-configuration $c$.

Let $C_{n,s,r}^{X,Y}$ be the set of all $(n,s,r)^{X,Y}$-configurations.  We claim that
\begin{equation*}
\left| C_{n,s,r}^{X,Y} \right| = \left| X_n^c \right| !{ \left| X_n^c \right|+r \choose r}
{n+1 \choose s-r} \prod\limits_{x \in X_n} (1+r+ \alpha_{X,n,x} + \beta_{Y,n,x}).
\end{equation*}
That is, we can construct  the $(n,s,r)^{X,Y}$-configurations as follows. 
First, we pick an order for the elements in $X_n^c$. This can be done in $\left| X_n^c \right|!$ ways. Next, we insert the $r$ $+$'s. This can be done in 
${ \left| X_n^c \right|+r \choose r}$ ways. Next, we insert the elements
of $X_n = \{ x_1 < x_2 < \cdots < x_{\left| X_n \right|} \}$ in increasing 
order. 
After placing $x_1, x_2, \ldots, x_{i-1}$, the next element $x_i$ can go 
\begin{itemize}
\item[$\bullet$] immediately before any of the $\beta_{Y,n,x_i}$ elements of $\{1, 2, \ldots, x_{i-1} \}$ that is not in 
$Y$, or 
\item[$\bullet$] immediately before any of the $\alpha_{X,n,x_i}$ elements of $\{x_i+1,x_i+2, \ldots, n\}$ that is not in $X$, or
\item[$\bullet$] immediately before any of the $r$ $+$'s, or
\item[$\bullet$] at the very end of the array.
\end{itemize}
Thus we can place the elements of $X_n$ in $\prod_{i=1}^{\left| X_n \right|} (1+r + \alpha_{X,n,x} + \beta_{Y,n,x})$ 
ways. Note that although $x_i$ might also be in $Y$, and might be placed immediately 
\emph{after} some other element of $X_n$, condition (ii) is not violated because the elements of $X_n$ 
are placed in increasing order. Finally, since each $-$ must occur either at the very start of the configuration 
or immediately following a number, we can place the $-$'s in ${n+1 \choose s-r}$ ways.

We define the weight $w(c)$ of an $(n,s,r)^{X,Y}$-configuration $c$ to be $(-1)^{s-r}$, i.e., $-1$ to the number 
of $-$'s of $c$. It then follows 
that the RHS of (\ref{eq:combXY1}) equals 
\begin{equation*}
\sum_{r=0}^s \sum_{c \in C_{n,s,r}^{X,Y}} w(c).
\end{equation*}
We now prove the theorem by exhibiting a sign-reversing involution $I$ on the set $C_{n,s}^{X,Y} = 
\bigsqcup\limits_{r=0}^s C_{n,s,r}^{X,Y}$, whose fixed points correspond to permutations $\sigma
\in S_n$ such that $des_{X,Y}(\sg) = s$. We say that a sign can be ``reversed''  if it can be changed from 
$+$ to $-$ or from $-$ to $+$ without violating conditions (i) and (ii). To apply $I$ to a configuration $c$, we scan from left
to right until we find the first sign that can be reversed. We then reverse that sign, and we let $I(c)$ be the resulting configuration.
If no signs can be reversed, we set $I(c) = c$.

In the example above, the first sign we encounter is the $+$ following $5$. This $+$ can be reversed, since 
$52$ is not an $X,Y$-descent. Thus $I(c)$ is the configuration shown below.
\begin{equation*}
I(c) = 5-2-+46+13-
\end{equation*}
It is easy to see that $I(I(c)) = c$ in this case, since applying $I$ again we change the $-$ following $5$ back to a $+$.

As another example, suppose $X=E, Y=O$, and $n=9$. Let $c$ be the following $(9,4,3)^{X,Y}$-configuration.
\begin{equation*}
c = 986+17-++4253
\end{equation*}
In this example we cannot reverse the $+$ following $6$, because $61$ is an $X,Y$-descent in the underlying permutation
$986174253$. Thus we move on to 
the $-$ following $7$. Changing this $-$ to a $+$ we get 
\begin{equation*}
I(c) = 986+17+++4253.
\end{equation*}
It is easy to see that $I(I(c))=c$ in this case, as well.

Conditions (i) and (ii) are clearly preserved by the very definition of $I$. It is also clear that $I$ is sign-reversing, 
since if $I(c) \not= c$, then $I(c)$ either has one more $-$ than $c$, or one fewer $-$ than $c$. 
To see that $I$ is in fact an involution, we note that the only signs that are \emph{not} reversible are single $+$'s 
occurring in the middle of an $X,Y$-descent pair, and $+$'s that immediately follows another sign. In either case, it is 
clear that a sign is reversible in a configuration $c$ if and only if the corresponding sign
is reversible in $I(c)$. Thus, if a sign is the first reversible sign in $c$, the corresponding sign in $I(c)$ must also be the first reversible sign in $I(c)$.
It follows that $I(I(c)) = c$ for all $c \in C_{n,s}^{X,Y}$. We therefore have 
\begin{equation*}\label{eq:combXY12}
\sum_{r=0}^s \sum_{c \in C_{n,s,r}^{X,Y}} w(c) = \sum_{r=0}^s \sum_{\small \begin{array}{cc} c 
\in C_{n,s,r}^{X,Y} \\ I(c) =c \end{array} \normalsize} w(c).
\end{equation*}
 
Now, consider the fixed points of $I$. Suppose that $I(c)=c$. Then $c$ clearly can have no $-$'s, and so $r=s$ and $w(c) = 1$.
It must also be the case that no $+$'s can be reversed. Thus each of the $s$ $+$'s must occur singly in the 
middle of an $X,Y$-descent pair. It follows that the underlying permutation has exactly $s$ $X,Y$-descents.

Finally, we should observe that if $\sigma = \sigma_1 \sigma_2 \cdots \sigma_n$ is a permutation 
with exactly $s$ $X,Y$-descents, then we can create a fixed point of 
$I$ simply by placing a $+$ in the middle of each $X,Y$-descent pair.
For example, if $X =\{2, 4, 6, 9\}$, $Y = \{1, 4, 7\}$, $n=9$, $s =2$, and 
$\sigma = 528941637$, then we have
\begin{equation*}
c = 5289+4+1637.
\end{equation*}
\end{proof}

Note that the right-hand side of (\ref{eq:combXY1}) makes sense for $s > \left| X_n \right|$, even though a permutation $\sigma \in 
S_n$ can have no more than $\left| X_n \right|$ $X,Y$-descents. This issue becomes important if one attempts to prove (\ref{eq:combXY1}) by induction, 
using the recursion (\ref{XYrec}). Although it is straightforward to show that the right-hand side of (\ref{eq:combXY1}) satisfies the 
recursion, to complete the proof one needs to show independently that the sum is zero when $s = \left| X_n \right|+1$. Our involution 
makes it clear that the sum is zero for \emph{all} $s > \left| X_n \right|$, since in such cases there will 
always be at least one $-$ or at least
$\left( \left| X_n \right|+1\right)$ $+$'s. If there are at least $\left( \left| X_n \right|+1\right)$ $+$'s, then there 
must be either an $X,Y$-descent $(c_i,c_{i+1})$ such that at least two $+$'s occur between $c_i$ and $c_{i+1}$, or
consecutive numbers $c_i$ and $c_{i+1}$ which do not form $X,Y$-descent, such that at least one $+$ appears between 
$c_i$ and $c_{i+1}$.
In each of these cases we can change a sign, and thus $I(c) \not= c$ for all $c \in C^{X,Y}_{n,s}$.

We can use the same involution on a related set of objects to prove an alternative formula for $P_{n,s}^{X,Y}$. 
\begin{theorem}\label{combXY2}
\begin{equation}\label{eq:combXY2}
P_{n,s}^{X,Y} = \left| X_n^c \right| !\sum_{r=0}^{\left| X_n \right|-s} (-1)^{\left| X_n \right|-s-r} \binom{\left| X_n^c \right|+r}{r}
\binom{n+1}{\left| X_n \right|-s-r}
\prod\limits_{x \in X_n} (r+\beta_{X,n,x}-\beta_{Y,n,x}).
\end{equation}
\end{theorem}

\begin{proof}
Let $X,Y,n$, and $s$ be given. For $r$ satisfying $0 \leq r \leq \left| X_n \right| -s$, an 
$\overline{(n,s,r)}^{X,Y}$-\emph{configuration} consists of an array of the numbers $1, 2, \ldots, n$, $r$ $+$'s, and $(\left| X_n \right|-s-r)$ $-$'s, 
satisfying the following three conditions.
\begin{enumerate}
\item[(i)] each $-$ is either at the very beginning of the array or immediately follows a number,
\item[(ii)] if $c_i \in X, 1 \leq i < n$, and $(c_i, c_{i+1})$ is not an $X,Y$-descent pair of the underlying 
permutation, then there must be at 
least one $+$ between $c_i$ and $c_{i+1}$, and
\item[(iii)] if $c_n \in X$, then at least one $+$ must occur to the right of $c_n$.
\end{enumerate}
Note that in an $\overline{(n,s,r)}^{X,Y}$-configuration, the number of $+$'s plus the number of $-$'s equals $\left| X_n \right|-s$.
As an example, if $X = \{2, 3, 6 \}$ and $Y = \{1,2,5\}$, then the following is a 
$\overline{(6,1,1)}^{X,Y}$-configuration.
\begin{equation*}
c = 213+6-54
\end{equation*}

Let $\overline{C}_{n,s,r}^{X,Y}$ be the set of all $\overline{(n,s,r)}^{X,Y}$-configurations.  Then we claim 
that 
\begin{equation*}
\left| \overline{C}_{n,s,r}^{X,Y} \right| = \left| X_n^c \right| !{ \left| X_n^c \right|+r \choose r}
{n+1 \choose \left| X_n \right| -s-r} \prod_{x \in X_n} (r+\beta_{X,n,x}-\beta_{Y,n,x}).
\end{equation*}
That is, we can construct  the $\overline{(n,s,r)}^{X,Y}$-configurations as follows. 
First, we pick an order for the elements in $X_n^c$. This can be done in $\left| X_n^c \right|!$ ways. Next, we insert the $r$ $+$'s. This can be done in 
${ \left| X_n^c \right|+r \choose r}$ ways. Next, we insert the elements
of $X_n = \{ x_1 < x_2 < \cdots < x_{\left| X_n \right|} \}$ in increasing 
order. First,
we can place $x_1$ in $r+\beta_{X,n,x_1}-\beta_{Y,n,x_1}$ ways, since $x_1$ can either go 
immediately before any of the $r$ $+$'s or 
immediately before any of the $x_1-1-\beta_{Y,n,x_1}=\beta_{X,n,x_1}-\beta_{Y,n,x_1}$ elements of $Y$ which are 
less than $x_1$. We note here that $\beta_{X,n,x_i} = x_i - i$ for all $i, 1 \leq i \leq \left| X_n\right|$. There are now two cases for placing $x_2$. \\
\ \\
{\bf Case 1.} $x_1$ was placed immediately in front of some element of $y \in Y$. In this case, 
$x_2$ cannot be placed immediately in front of $y$, 
since otherwise we would violate condition (ii). $x_2$ can be placed before  
any $+$ or immediately in front of any element of $Y$ which is less than $x_2$, except 
$y$. Hence, $x_2$ can be placed in
\begin{eqnarray*}
r+x_2-1-\beta_{Y,n,x_2}-1 & =& r+x_2-2-\beta_{Y,n,x_2}\\
 & = & r+\beta_{X,n,x_2}-\beta_{Y,n,x_2}
\end{eqnarray*}
ways.\\
\ \\
{\bf Case 2.} $x_1$ was placed immediately before a $+$. In this case, $x_2$ cannot be placed immediately before 
the same $+$, since again we would violate condition (ii).  $x_2$ can be placed immediately before any of the 
other $+$'s or immediately before any element of $Y$ which is less than $x_2$. Hence  $x_2$ can be 
placed in
\begin{eqnarray*}
r-1+x_2-1-\beta_{Y,n,x_2}& =& r+x_2-2-\beta_{Y,n,x_2}\\
 & = & r+\beta_{X,n,x_2}-\beta_{Y,n,x_2}
\end{eqnarray*}
ways.\\
\ \\
In general, having placed 
$x_1, x_2, \ldots, x_{i-1}$, we cannot place $x_i$ immediately before 
some $y \in Y, y < x_i$, which earlier had an element of $\{x_1, x_2, \ldots, x_{i-1}\}$ 
placed before it. Similarly, we cannot place $x_i$ immediately before any $+$ which earlier had 
an element of $\{x_1, x_2, \ldots, x_{i-1}\}$ placed before it. It then follows that there are 
\begin{eqnarray*}
r+x_i-1-\beta_{Y,n,x_i}-(i-1) & = & r+x_i-i-\beta_{Y,n,x_i} \\
 & = & r+\beta_{X,n,x_i}-\beta_{Y,n,x_i}
\end{eqnarray*}
ways to place $x_i$. Thus, there are total of $\prod_{i=1}^{\left| X_n\right|} 
(r+\beta_{X,n,x_i}-\beta_{Y,n,x_i})$ ways to place $x_1, x_2, \ldots, x_{\left| X_n\right|}$, 
given our placement of the elements of $X_n^c$.
Finally, we can place the $-$'s in ${n+1 \choose \left| X_n\right|-s-r}$ ways. 

We define the weight $w(c)$ of an $\overline{(n,s,r)}^{X,Y}$-configuration $c$ to be $(-1)^{\left| 
X_n\right|-s-r}$, 
i.e., $-1$ to the number of $-$'s of $c$. It then follows that the RHS of (\ref{eq:combXY2}) equals 
\begin{equation*}
\sum_{r=0}^{\left| X_n \right|-s} \sum_{c \in \overline{C}_{n,s,r}^{X,Y}} w(c),
\end{equation*}
We now prove the theorem by exhibiting a sign-reversing involution $I$ on the set $\overline{C}_{n,s}^{X,Y} = \bigsqcup\limits_{r=0}^{\left| X_n\right|-s} 
\overline{C}_{n,s,r}^{X,Y}$ whose fixed points correspond to permutations $\sigma \in S_n$ such that $des_{X,Y}(\sg) = s$.
We define $I$ exactly as in the proof of Theorem \ref{combXY1}. That is, we scan from left to right and reverse the 
first sign that we can reverse without violating conditions (i)-(iii).

In the example above, we cannot reverse the $+$ following $3$ without violating condition (ii), since $3 \in X$ and $36$ is not an $X,Y$-descent.
Thus, we reverse the $-$ following $6$ to get
\begin{equation*}
I(c) = 213+6+54.
\end{equation*}

We argue as in Theorem \ref{combXY1} that $I$ is a sign-reversing involution, so that
\begin{equation*}
\sum_{r=0}^s \sum_{c \in \overline{C}_{n,s,r}^{X,Y}} w(c) = \sum_{r=0}^s \sum_{\small \begin{array}{cc} 
c \in \overline{C}_{n,s,r}^{X,Y} \\ I(c) =c\end{array} \normalsize} w(c).
\end{equation*}
 
Now, consider a fixed point $c$ of $I$. As in the proof of Theorem \ref{combXY1}, $c$ can have no $-$'s, and 
thus $r=\left| X_n\right|-s$ and $w(c) = 1$. No 
string of multiple $+$'s can occur, since the first $+$ in such a string could be reversed. Thus, each of the $\left( 
\left| X_n\right|-s \right)$ $+$'s appears 
singly, and must either
\begin{itemize}
\item[$\bullet$] immediately follow some $c_i \in X, 1 \leq i < n$, such that $(c_i,c_{i+1})$ is not an 
$X,Y$-descent pair of the underlying permutation, or
\item[$\bullet$] immediately follow $c_n\in X$.
\end{itemize}
Thus $\left| X_n\right|-s$ elements of $X_n$ immediately precede a $+$ that cannot be reversed, and are thus not 
the tops of $X,Y$-descent pairs. 
It follows that each of the remaining $s$ elements of $X_n$ do not immediately precede a $+$, 
and as such each must be the top of an $X,Y$-descent pair,. Thus the underlying 
permutation $c_1 c_2 \cdots c_n$ has exactly $s$ $X,Y$-descents. 

Again, we observe that if $\sigma_1 \sigma_2 \cdots \sigma_n$ is a permutation 
with exactly $s$ $X,Y$-descents, then we can create a fixed point of 
$I$ by inserting a $+$
after every element of $X_n$ that is not 
the top of an $X,Y$-descent pair. For example, 
if $X =\{2, 3, 4, 6, 8,9\}$, $Y = \{1,2,3,5\}$, $n=9$, $s =4$, and 
$\sigma =958621437$, then the corresponding configuration would be
\begin{equation*}
958+62143+7.
\end{equation*}
\end{proof}

We note that the 
quantity $r + \beta_{X,n,x_i}-\beta_{Y,n,x_i}$ may be zero, or even negative. 
For example, let $X = \{2, 3, 4 \}$, $Y = \{1, 3, 5\}$, $n 
= 6$, and $s=3$. 
In constructing a $\overline{(6,3,0)}^{X,Y}$ configuration, we  
start with an ordering of $X_6^c = \{1, 5, 6 \}$, such as
\begin{equation*}
516.
\end{equation*}
Since $\left| X_6 \right| - s = 3-3 = 0$, there are no $+$'s to place, and the next step is to place the elements of 
$X_6$ in increasing order.
There is one place to put $x_1 = 2$, namely, immediately before the $1$. This corresponds to the 
fact that 
$\beta_{X,6,2}-\beta_{Y,6,2} = (2-1)-0 = 1$. We then have the array
\begin{equation*}
5216.
\end{equation*}
Notice that now there is no place to put $x_2 = 3$. There are no $+$'s, and $3$ cannot be placed in front of $1$ (the only element 
of $Y$ smaller than $3$) without violating condition (ii). This corresponds to the fact that 
$\beta_{X,6,3}-\beta_{Y,6,3} = (3-2)-1 = 0$. So $P^{X,Y}_{6,3} = 0$, which we can also see by inspection: the only potential 
$X,Y$-descents are 
$21, 31, 41$, and $43$, and no permutation can contain more than one of $21, 31$, and $41$. Note, finally, that the quantity $r + 
\beta_{X,n,x_i}-\beta_{Y,n,x_i}$ is non-negative for $i=1$, since $\beta_{Y,n,x_1} \leq \beta_{X,n,x_1} = x_1-1$, 
and that the difference 
\begin{eqnarray*}
&& (r + \beta_{X,n,x_i}-\beta_{Y,n,x_i}) - (r + \beta_{X,n,x_{i+1}}-\beta_{Y,n,x_{i+1}}) \\
&&
 = (r + x_i-i -\beta_{Y,n,x_i}) - (r + x_{i+1}-(i+1)-\beta_{Y,n,x_{i+1}}) \\
 && = 1 + x_i - \beta_{Y,n,x_i}-(x_{i+1}-\beta_{Y,n,x_{i+1}} )
\end{eqnarray*}
is at 
most $1$, since the sequence $\{ x_i - \beta_{Y,n,x_i} \}_{i=1}^{\left| X_n\right|}$ is nondecreasing. Thus in a situation 
in which $r + \beta_{X,n,x_i}-\beta_{Y,n,x_i}$ is negative, we must have 
$r + \beta_{X,n,x_j}-\beta_{Y,n,x_j} = 0$ for some $j < i$.

Taking $Y = \mathbb{N}$ in (\ref{eq:combXY1}) and (\ref{eq:combXY2}) gives the following two corollaries.

\begin{corollary}\label{comb1}
\begin{equation}\label{eq:comb1}
P_{n,s}^X = \left| X_n^c\right| !\sum_{r=0}^s (-1)^{s-r} \binom{\left| X_n^c\right|+r}{r}
\binom{n+1}{s-r}
\prod\limits_{x \in X_n} (1+r + \alpha_{X,n,x})
\end{equation}
\end{corollary}

\begin{corollary}\label{comb2}
\begin{equation}\label{eq:comb2}
P_{n,s}^X = \left| X_n^c\right| !\sum\limits_{r=0}^{\left| X_n\right|-s} (-1)^{\left| X_n\right|-s-r} {\left| X_n^c\right|+r \choose r}
{n+1 \choose \left| X_n\right|-s-r}
\prod\limits_{x \in X} (r+\beta_{X,n,x})
\end{equation}
\end{corollary}

\begin{remark}
Using a rook-placement interpretation outlined in Section 4 of this paper, (\ref{eq:comb1}) and (\ref{eq:comb2}) can also 
be obtained from formulas of Haglund \cite{Haglund} for hit numbers of Ferrers boards. 
\end{remark}

One can use (\ref{eq:combXY1}) and (\ref{eq:combXY2}) to obtain similar formulas for the coefficients of the 
polynomials $Q_n^X(x) = P_n^{\mathbb{N},X}(x)$.
Formulas for the coefficients of $Q_n^X(x)$ can also be derived directly from (\ref{eq:comb1}) and (\ref{eq:comb2}) by using the 
following result. 

\begin{theorem}
Given a subset $X \subseteq \mathbb{N}$ and a permutation $\sigma \in S_n$, let $X^*$ be the 
subset of $[n]$ satisfying $i \in X^* \iff n+1-i \in X_n$. Then
\begin{equation*}
Q_{n,s}^X = P_{n,s}^{X^*}. 
\end{equation*}
\end{theorem}
\begin{proof}
Given a permutation $\sigma = \sigma_1 \sigma _2 \ldots \sigma _n$, the \emph{complement} of $\sigma$ is 
$\sigma^c = (n+1-\sigma_1)(n+1- \sigma _2) \ldots (n+1-\sigma _n)$. The \emph{reverse} of $\sigma$ is $\sigma^r = 
\sigma_n \ldots \sigma_2 \sigma_1$. The operations of ``complement'' and ``reverse'' are both clearly 
invertible. Now suppose $\sigma$ has a descent pair $(i,j)$. Then $\sigma^c$ has an ascent pair $(n+1-i,n+1-j)$, and 
so $\left(\sigma^c \right)^r$ has a descent pair $(n+1-j,n+1-i)$. Thus 
$\overrightarrow{des}_X(\sigma) = \overleftarrow{des}_{X^*}\left( \left( \sigma^c 
\right)^r \right)$, which implies $Q_{n,s}^X = P_{n,s}^{X^*}$.
\end{proof}

For example, if  $X=k\mathbb{N}$ and $n = km+j$ for some $j, 0 \leq j \leq k-1$, then $X_n 
= \{k,2k, \ldots km\} \subseteq [n]$ and $X^* = \{1+j,1+j+k, 
\ldots ,1+j+k(m-1)\}$. Thus $\left| (X^*)_n^c \right| = (k-1)m+j$, and, for $i=0, \ldots, m-1$, we have
\begin{eqnarray*}
\alpha_{X^*,n,1+j+ik} &=& km+j-(1+j+ik)-(m-1-i) = (k-1)(m-i), \mbox{ and}\\
\beta_{X^*,n,1+j+ik} &=& j+ (k-1)i.
\end{eqnarray*}
It follows from Corollaries \ref{comb1} and \ref{comb2} that 
\begin{eqnarray*}
&&Q_{km+j,s}^{k\mathbb{N}} = P_{km+j,s}^{\{1+j,1+j+k, \ldots ,1+j+k(m-1)\}} = \\
&& \hspace{0.6cm} ((k-1)m+j)!\sum_{r=0}^s (-1)^{s-r}\binom{(k-1)m+j +r}{r} \binom{km+j+1}{s-r} 
\prod_{i=1}^m (1+r+(k-1)i), \mbox{ and}\\
&&Q_{km+j,s}^{k\mathbb{N}} = \\
&& \hspace{0.6cm} ((k-1)m+j)!\sum_{r=0}^{m-s} (-1)^{m-s -r}\binom{(k-1)m+j +r}{r} \binom{km+j+1}{m-s-r} 
\prod_{i=0}^{m-1} (r+j+(k-1)i).
\end{eqnarray*} 
These two formulas for $Q_{km+j,s}^{k\mathbb{N}}$ were first proved by 
Kitaev and Remmel in \cite{KitaevRemmel2} using the special case of the recursion (\ref{XYrec}) with
$X =\mathbb{N}$ and $Y = k \mathbb{N}$.

\section{Applications}
\label{cases-section}

\begin{corollary}
Let $X = \mathbb{N}$, so that $\overleftarrow{Des}_X(\sigma) = Des(\sigma)$. Then
\begin{equation*}
P_{n,s}^X = \sum\limits_{r=0}^s (-1)^{s-r}{n+1 \choose s-r} (1+r)^n,
\end{equation*}
which is a well-known formula for the Eulerian numbers (see, e.g.,  \cite{Comtet}, pp. 240--246).
\end{corollary}

In the following results we employ the notation of hypergeometric series. For $a \in \mathbb{R}$ and $n \in \mathbb{N}$, let
$(a)_n = a(a+1)(a+2)\cdots (a+n-1)$. Let $(a)_0 = 1$. Define
\begin{equation*}
_{m+1}F_m \left[ \begin{array}{rrrrr} a_0, 
& a_1, & a_2, & \ldots, & a_m \\ & b_1, & b_2, & \ldots, & b_m \end{array} \right]
:= \sum\limits_{r=0}^\infty \frac{(a_0)_r (a_1)_r (a_2)_r \cdots (a_m)_r}{r! (b_1)_r (b_2)_r \cdots (b_m)_r}.
\end{equation*}
Since $(-n)_r = 0$ for all $r > n$, a hypergeometric series may be undefined if a parameter in the 
denominator is a negative integer. In our applications, \emph{all} of the parameters are negative integers. 
However, in each case the largest (least negative) parameter occurs in the numerator, hence the series terminates 
in a well-defined, finite sum.

\begin{corollary}\label{ap1}
Let $X = 2\mathbb{N}$. Then
\begin{equation*}
P_{2n,s}^X = (n!)^2{n \choose s}^2,
\end{equation*}
which was originally derived by Kitaev and Remmel \cite{KitaevRemmel1} using the special case of the 
recursion (\ref{XYrec}) with
$X = 2 \mathbb{N}$ and $Y = \mathbb{N}$.
\end{corollary}
\begin{proof}
By Corollary \ref{comb1}, we have
\begin{eqnarray*}
P_{2n,s}^X & = & n!\sum\limits_{r=0}^s (-1)^{s-r} {n+r 
\choose r} {2n+1 \choose s-r} \prod\limits_{i=1}^n(i+r) \\
 & = & (s+1)^2_n \mbox{ }_3F_2 \left[ \begin{array}{rrr} -s, 
& -s, & -(2n+1) \\ & -(n+s), & -(n+s) \end{array} \right] \\
 & = & (s+1)^2_n \frac{(n+1-s)_s (n+1-s)_s}{(n+1)_s (n+1)_s} \\
 & = & (n!)^2{n \choose s}^2,
\end{eqnarray*}
where in the third step we use the Pfaff-Saalsch\"{u}tz $\mbox{}_3F_2$ summation 
formula (see \cite{GasperRahman})
\begin{equation*}
\mbox{ }_3F_2 \left[ \begin{array}{rrr} -n, & a, & b \\ & c, 
& a+b-c-n+1 \end{array} \right] = \frac{(c-a)_n(c-b)_n}{(c)_n(c-a-b)_n}.
\end{equation*}
\end{proof}

\begin{remark}\label{ap2}
More generally, if $X = \{u+2, u+4, u+6, \ldots, u+2m \}$, a similar computation gives
\begin{equation}\label{eq:ap2}
P_{2m+u+v,s}^X = {m \choose s} {m+u+v \choose v+s}(m+u)!(m+v)!.
\end{equation}
\end{remark}

Combining (\ref{eq:comb1}) and (\ref{eq:comb2}) for various sets $X$, we get interesting identities, such as the 
following result of 
Kitaev and Remmel (Theorem 2 of \cite{KitaevRemmel2}).
\begin{corollary}
Let $X = k\mathbb{N}$. Then for each $j, 0 \leq j \leq k-1$, we have
\small
\begin{eqnarray*}
&&P_{kn+j,s}^X = \\
&& \hspace{0.6cm} ((k-1)n+j)! \sum\limits_{r=0}^s (-1)^{s-r}{(k-1)n+j + r \choose r}{kn+j+1 \choose s-r} \prod\limits_{i=0}^{n-1}
\left(1+r+j+(k-1)i\right), \mbox{ and}\\
&&P_{kn+j,s}^X = \\
&& \hspace{0.6cm} ((k-1)n+j)! \sum\limits_{r=0}^{n-s} (-1)^{n-s-r}{(k-1)n+j + r \choose r}{kn+j+1 \choose n-s-r} \prod\limits_{i=1}^n
\left(r+(k-1)i\right).
\end{eqnarray*}
\normalsize
\end{corollary}

For some sets $X$ the right-hand sides of (\ref{eq:comb1}) and (\ref{eq:comb2})
can be rewritten in terms of hypergeometric series; hence we obtain combinatorial
proofs of identities such as the following.
\begin{corollary}
Let $k$ and $m$ be positive integers, and let $s$ be a non-negative integer. Then
\begin{eqnarray*}
&&(s+1)_m^{k+1} \mbox{ }_{k+2}F_{k+1} \left[ \begin{array}{cccc} 
-((k+1)m+1), & -s, & \ldots, & -s \\
 & -(m+s), & \ldots, & -(m+s) \end{array} \right] = \\
&& \hspace{0.2cm} (km+1-s)_m^{k+1} \mbox{ }_{k+2}F_{k+1} \left[ 
\begin{array}{cccc}
-((k+1)m+1), & -(km-s), & \ldots, & -(km-s) \\
& -((k+1)m-s), & \ldots, & -((k+1)m-s) \end{array} \right].
\end{eqnarray*}
\end{corollary}

\begin{proof}
Let $X = \{ i : i \not= 1 \mod (k+1)\}$, and use (\ref{eq:comb1}) and (\ref{eq:comb2}) to compute $P_{(k+1)m,s}^X.$ 
\end{proof}

This identity is a special case of an integral form of a 
transformation of Karlsson-Minton type hypergeometric series due to 
Gasper \cite{Gasper}:

\begin{equation*}
\hspace{-2cm}\mbox{ }_{k+2}F_{k+1} \left[ \begin{array}{ccccc}
 w, & x, & b_1+d_1, & \ldots, & b_k + d_k \\
 & x+c+1, & b_1, & \ldots, & b_k \end{array} \right] = 
\end{equation*}
\begin{equation*}
\hspace{-3cm}\frac{\Gamma(1+x+c)\Gamma(1-w)}{\Gamma(1+x-w)\Gamma(c+1)}\prod\limits_{i=1}^k 
\frac{(b_i-x)_{d_i}}{(b_i)_{d_i}}
\end{equation*}
\begin{equation*}
\hspace{2cm}\times \mbox{ }_{k+2}F_{k+1} \left[ \begin{array}{ccccc}
 -c, & x, & 1+x-b_1, & \ldots, & 1+x-b_k \\
 & x+1-w, & 1+x-b_1-d_1, & \ldots, & 1+x-b_k-d_k \end{array} \right] 
\end{equation*}

The following identity, while still a special case of Gasper's transformation, is more general.

\begin{corollary}
Let $u = (u_1, \ldots, u_k)$ be a weakly increasing array of non-negative integers, and $v = (v_1, \ldots, v_k)$ an array of positive integers.
Then for $n \geq \sum\limits_{i=1}^k v_i + max \{u_i + v_i : 1 \leq i \leq k \} - u_1$, we have
\begin{eqnarray*}
& & (s+1)_a (s+u_1+1)_{v_1} \cdots (s+u_k+1)_{v_k} \\
& & \hspace{1cm} \times \mbox{ }_{k+2}F_{k+1} \left[ \begin{array}{ccccc} 
-(n+1), & -s, & -(s+u_1) & \ldots, & -(s+u_k) \\
 & -(s+a) & -(s+u_1+v_1) & \ldots, & -(s+u_k+v_k) \end{array} \right] =
\end{eqnarray*} 
\begin{eqnarray}\label{bal}
 &&  (-1)^n(-n+s)_a  (-n+s+u_1)_{v_1} \cdots (-n+s+u_k)_{v_k} \nonumber \\
 & & \hspace{1cm} \times \mbox{ }_{k+2}F_{k+1} 
\left[ 
\begin{array}{ccccc}
-(n+1), & -n + s + a & -n + s + u_1 + v_1 & \ldots, & -n+s+u_k+v_k \\
 & -n+s & -n +s + u_1 & \ldots, & -n+s+u_k \end{array} \right], \nonumber \\
 & & 
\end{eqnarray}
where $a = n-\sum\limits_{i=1}^k v_i$.
\end{corollary}

\begin{remark}
Note that both hypergeometric series in (\ref{bal}) are \emph{balanced}, i.e., the sum of the parameters in the 
top row is one less than the sum of the parameters in the bottom row. Balanced hypergeometric series are a 
particularly well-behaved class of hypergeometric series for which several summation and transformation results exist.
\end{remark}

\begin{proof}
For each $i, 1 \leq i \leq k$, let
\begin{equation*}
f(i) = \left| \{ j : u_j \leq i \leq u_j + v_j -1 \} \right| .
\end{equation*}
Define $M = \mbox{max} \{ m : f(m) > 0 \}$ and set $b = a+1-M-u_1$.
Let $X$ be the subset of $\mathbb{N}$ defined by the binary sequence
\begin{equation*}
\tau = \underbrace{0 \ldots 0}_b \underbrace{1\ldots 1}_{f(M)} 0 \underbrace{1\ldots 1}_{f(M-1)} 0 \underbrace{1\ldots 1}_{f(M-2)} 0 
\ldots 0 \underbrace{1\ldots 1}_{f(u_1+1)} 0 
\underbrace{1\ldots 1}_{f(u_1)} \underbrace{0\ldots 0}_{u_1}.
\end{equation*}
That is, let $i \in X$ if and only if $\tau_i = 1$.

We prove the identity by showing that both sides are equal to $P_{n,s}^X$. Applying (\ref{eq:comb1}) gives
\begin{eqnarray*}
P_{n,s}^X & = & a!\sum\limits_{r=0}^s (-1)^{s-r} {a+r \choose r}
{n+1 \choose s-r} \prod\limits_{i=u_1}^M (i+1+r)^{f(i)} \\
 &  = & a!\sum\limits_{r=0}^s (-1)^{s-r} {a+r \choose r}
{n+1 \choose s-r} \prod\limits_{i=1}^k (u_i+1+r)_{v_i} \\
 & = & (s+1)_a (s+u_1+1)_{v_1} \cdots (s+u_k+1)_{v_k} \\
  & & \hspace{1cm} \times \mbox{ }_{k+2}F_{k+1} \left[ \begin{array}{ccccc} 
-(n+1), & -s, & -(s + u_1) & \ldots, & -(s + u_k) \\
 & -(s+a) & -(s+u_1+v_1) & \ldots, & -(s+u_k+v_k) \end{array} \right].
\end{eqnarray*}
The key observation here is that $\prod\limits_{i=u_1}^M (i+1+r)^{f(i)} = \prod\limits_{i=1}^k (u_i+1+r)_{v_i}$, since
$i+1+r$ occurs in $\prod\limits_{j=1}^k (u_j+1+r)_{v_j}$ exactly 
$f(i) = \left| \{ j : u_j \leq i \leq u_j + v_j -1 \} \right|$ times.

On the other hand, applying (\ref{eq:comb2}) gives
\begin{eqnarray*}
P_{n,s}^X & = & a!\sum\limits_{r=0}^{n-a-s} (-1)^{n-a-s-r} {a+r \choose r}
{n+1 \choose n-a-s-r} \prod\limits_{i=u_1}^M (a - i +r)^{f(i)} \\
 & = & a!\sum\limits_{r=0}^{n-a-s} (-1)^{n-a-s-r} {a+r \choose r}
{n+1 \choose n-a-s-r} \prod\limits_{i=1}^k (a +1 - u_i -v_i +r)_{v_i} \\
 & = & (-1)^n(-n+s)_a  (-n+s+u_1)_{v_1} \cdots (-n+s+u_k)_{v_k}  \\
  &  & \hspace{1cm} \times \mbox{ }_{k+2}F_{k+1} 
\left[ 
\begin{array}{ccccc}
-(n+1), & -n + s + a & -n + s + u_1 + v_1 & \ldots, & -n+s+u_k+v_k \\
 & -n+s & -n +s + u_1 & \ldots, & -n+s+u_k \end{array} \right].
\end{eqnarray*}
\end{proof}

\begin{example}
Let $u = (0,1,1,5), v = (2,3,1,2)$, and $n = 16$. If we represent $u$ and $v$ with four rows of $\bullet$'s, the 
$i^{\mbox{th}}$ row
starting at $u_i$ and having length $v_i$, then $f(i) = \left| \{ j : u_j \leq i \leq u_j + v_j -1 \} \right|$ is the
number of $\bullet$'s in the $i^{\mbox{th}}$ column. In this example we have $f(0) = 1, f(1) = 3, f(2) = f(3) = 1, f(4) 
= 0$, and $f(5) = f(6) = 1$, as shown below.
\begin{equation*}
\begin{array}{ccccccc}
 & & & & & \bullet & \bullet \\
 & \bullet & & & & & \\
 & \bullet & \bullet & \bullet & & & \\
\bullet & \bullet & & & & & \\
0 & 1 & 2 & 3 & 4 & 5 & 6
\end{array}
\end{equation*}
The corresponding binary sequence is
\begin{equation*}
\tau = 0010100101011101,
\end{equation*}
and the corresponding set is $X = \{3, 5, 8, 10, 12, 13, 14, 16 \}$. Using (\ref{eq:comb1}) and 
(\ref{eq:comb2}) to compute
$P_{n,s}^X$ gives
\begin{eqnarray*}
& & (s+1)_{8} (s+1)_2 (s+2)_3 (s+2)_1 (s+6)_2 \\
 & & \hspace{1cm} \times \mbox{ }_{6}F_{5} \left[ \begin{array}{cccccc} 
-17, & -s, & -s, & -(s+1), & -(s+1), & -(s+5) \\
 & -(s+8), & -(s+2), & -(s+4), & -(s+2), & -(s+7)  \end{array} \right] 
\end{eqnarray*} 
\begin{eqnarray*} 
 & = &  (-1)^{16}(-16+s)_{8}  (-16+s)_2 (-15+s)_3 (-15+s)_1 (-11+s)_2 \\
 & & \hspace{1cm} \times \mbox{ }_{6}F_{5} \left[ 
\begin{array}{cccccc}
-17, & -8 + s, & -14 + s, & -12+s, & -14+s, & -9+s \\
 & -16+s & -16 +s, & -15+s, & -15+s, & -11+s \end{array} \right].
\end{eqnarray*}
\end{example}


\section{Connections with Rook Theory}
\label{rook-section}

A \emph{board} is a finite subset of an infinite grid of unit squares. The ``rook number'' $r_k(B)$ of a board $B$ is defined to be
the number of ways to place $k$ non-attacking rooks on $B$. Two boards $B$ and $B'$ are \emph{rook-equivalent} if $r_k(B) = r_k(B')$ for 
all $k$.
For a board $B$ contained inside the $n \times n$ board, the
``hit number'' $h_k(B)$ is defined to be the number of ways to place
$n$ non-attacking rooks on the $n \times n$ board so that exactly $k$ rooks lie on $B$.
In what follows we will focus on hit numbers rather than rook numbers. Kaplansky and Riordan \cite{KaplanskyRiordan} showed
that rook-equivalent boards have the same hit numbers.

The key to the connection between rook placements and descents of permutations
is Foata's First Transformation \cite{FoataSchutzenberger1}, a bijection $\Phi: S_n 
\longrightarrow S_n$
which exchanges excedences and descents. An \emph{excedence} of $\sigma = \sigma_1 \sigma_2 \ldots \sigma_n$ is an entry
$\sigma_i$ satisfying $\sigma_i > i$. Foata's transformation can most easily be explained with an example.

\begin{example}
Let $\omega = 61437258 = \left( \begin{array}{cccccccc} 1 & 2 & 3 & 4 & 5 & 6 & 7 & 8\\
6 & 1 & 4 & 3 & 7 & 2 & 5 & 8 \end{array} \right)$. This permutation has three excedences: $\begin{array}{c} 1 \\ 6 \end{array},
\begin{array}{c} 3 \\ 4 \end{array}$, and $\begin{array}{c} 5 \\ 7 \end{array}$. The first step in Foata's transformation is to write $\omega$
in cycle form: $(162)(34)(57)(8)$. Next, write each cycle with largest element last, and order the cycles by increasing largest element:
$(34)(216)(57)(8)$. Finally, to compute $\Phi(\omega)$, reverse each cycle and erase the parentheses: $\Phi(\omega) = 43612758$. 
In this example the descents of $\Phi(\omega)$ are $43, 61$, and $75$. In general, it is not hard to see that $(i,j)$ is a descent pair of $\Phi(\omega)$ if and only if 
$\begin{array}{c} j \\ i \end{array}$ is an excedence of $\omega$.
To go backwards, given $\sigma = 43612758$, cut before each left-to-right maxima: $43|612|75|8$, then reverse each
block to get the cycles of $\Phi^{-1}(\sigma)$: $(34)(216)(57)(8)$.
\end{example}

Foata's transformation is key to this section because rook placements provide a convenient way of 
tracking the excedences of a permutation. As the following example illustrates, given any subset $U \subseteq
\{(i,j) : 1 \leq j < i \leq n \}$ of potential excedences we can construct
a board $B^U_n$ inside the $n \times n$ board so that the number of $\sigma \in S_n$ with exactly $s$ $U$-excedences, and hence 
the number of $\sigma \in S_n$ with exactly
$s$ $U$-descents, is $h_{s}(B^U_n)$.

\begin{example}\label{XYZrook}
Suppose we wish to count descents $\sigma_i > \sigma_{i+1}$ satisfying $\sigma_i \in E, \sigma_{i+1} \in O$, and $\sigma_i - \sigma_{i+1} \in \{1, 3\}$
(this is an instance of counting what we have called ``$X,Y,Z$-descents''). For $n=8$, the board 
$B^U_8$ consists of the squares $(i,j) \in [8] \times [8]$ 
such that $i \in E, j  \in O$, and  $i-j=1$ or $3$.  We have pictured this board as the shaded squares 
in Figure \ref{fig:fig16}.
\fig{.75}{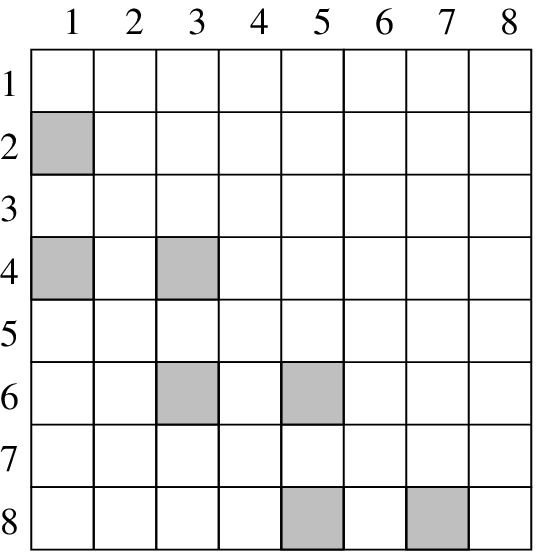}{The board $B^U_8$.}
Now consider the placement, shown in Figure \ref{fig:fig17}, of eight non-attacking rooks (marked by $X$'s) on the $8 
\times 8$ board so that two rooks lie on $B^U_8$.
\fig{.75}{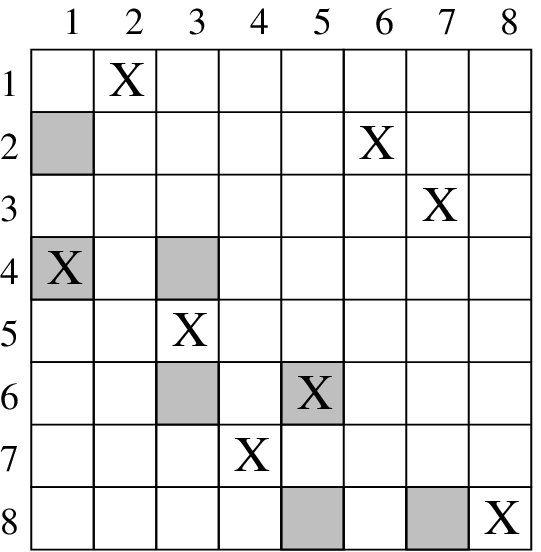}{A placement of rooks on the $8\times 8$ board.}
This placement corresponds to the permutation $\omega = \left( \begin{array}{cccccccc} 1 & 2 & 3 & 4 & 5 & 6 & 7 & 8\\
4 & 1 & 5 & 7 & 6 & 2 & 3 & 8 \end{array} \right)$, with the rooks placed on $B^U_8$ corresponding to the excedences 
$\begin{array}{c} 1 \\ 4 \end{array}$ and $\begin{array}{c} 5 \\ 6 \end{array}$. We now employ Foata's 
transformation to get the permutation $\sigma = \Phi(\omega) = 74126538$ with exactly two $U$-descents: $41$ and $65$.
\end{example}

One important class of boards is the class of \emph{Ferrers} boards, that is, boards of partition shape. 
Ferrers boards are usually drawn right justified, as in Figure \ref{fig:fig1}.
\fig{.75}{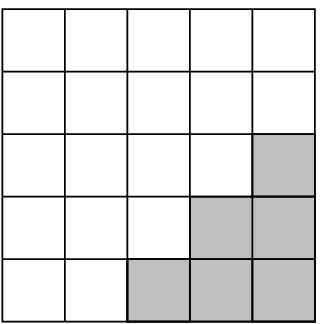}{A Ferrers board.}
For $X,Y \subseteq [n]$, the board $B_n^{X,Y}$ corresponding to the potential $X,Y$-descents of permutations $\sigma 
\in S_n$ is trivially 
rook-equivalent to 
a Ferrers board; we need only shift all the non-empty rows and columns to the bottom-left of the square and take the mirror image.
For example, if $X = \{2, 3, 5\}, Y = \{1,2, 4\}$, and $n=5$, the board $B_5^{X,Y}$ is shown in Figure 
\ref{fig:fig2}; this board is trivially rook-equivalent to the Ferrers board shown in Figure \ref{fig:fig1}.
\fig{.75}{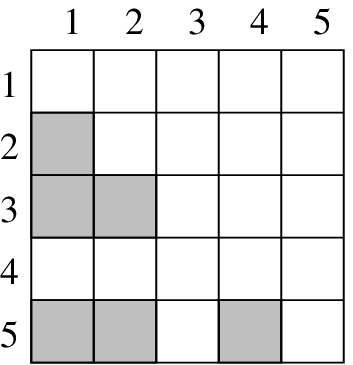}{The board $B^{X,Y}_5$.}
Therefore, in certain cases we can make use of results for Ferrers boards in computing the numbers $P_{n,s}^{X,Y}$. 
As one example, in \cite{Haglund}, Haglund gives several formulas involving the hit numbers of Ferrers
boards, one of which can be specialized to obtain (\ref{eq:comb1}). 

As another example, we 
can use the rook interpretation of our problem to give a purely combinatorial proof of Corollary \ref{ap1}. For 
instance, the board $B^{E,\mathbb{N}}_8$ corresponding to even descents of permutations $\sigma \in S_8$ is shown in 
Figure \ref{fig:fig3}.
\fig{.75}{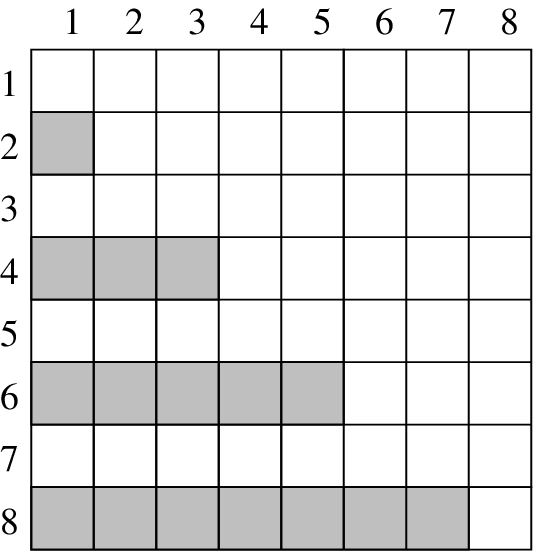}{The board $B^{E,\mathbb{N}}_8$.}
The corresponding Ferrers board is shown in Figure \ref{fig:fig4}.
\fig{.75}{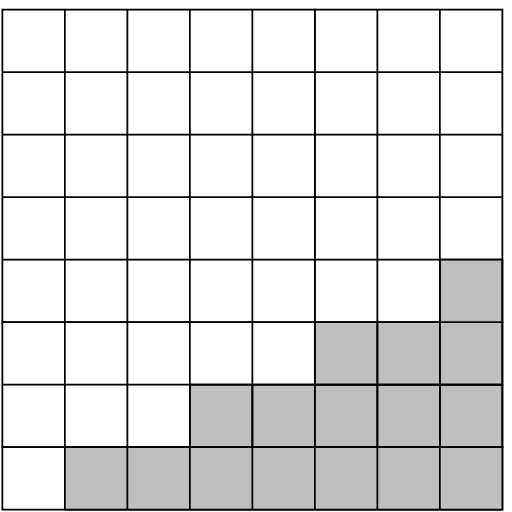}{The Ferrers board $B$ corresponding to $B^{E,\mathbb{N}}_8$.}

In what follows we use the notation of \cite{FoataSchutzenberger2} and \cite{GoldmanJoichiWhite}.
Let the column heights of a Ferrers board $B$ inside an $n \times n$ square be given by the ``height vector'' $h(B) = (h_1, h_2, \ldots, h_n)$. Define 
the ``structure vector'' $s(B) = (s_1, s_2, \ldots, s_n)$, where $s_i = h_i - (i-1), 1 \leq i \leq n$. Here it is standard practice to insist that $n$ be 
large enough so that none of the entries of the structure vector is positive. This can always be done, for example, by taking
$n$ greater than the number of squares of $B$. In our applications $n$ is already fixed; however, 
the entries of the structure vector are still non-positive because boards corresponding to descents necessarily
lie strictly below the main diagonal.
In \cite{FoataSchutzenberger2}, Foata and Sch\"{u}tzengerger showed that two Ferrers boards
$B$ and $B'$ are rook-equivalent if and only if the entries of $s(B)$ and $s(B')$ are equal as multisets.

In our example, we first compute $s(B)$:
\begin{center}
\begin{tabular}{rrrrrrrrr}
& $0$ & $1$ & $1$ & $2$ & $2$ & $3$ & $3$ & $4$ \\
$-$ & $0$ & $1$ & $2$ & $3$ & $4$ & $5$ & $6$ & $7$ \\
\hline
 & $0$ & $0$ & $-1$ & $-1$ & $-2$ & $-2$ & $-3$ & $-3$
\end{tabular}
\end{center}
Thus $B$ is rook-equivalent to the board $B^\prime$ with structure vector
$s(B') = (0,-1,-2,-3,0,-1,-2,-3)$. To identify $B^\prime$ we next compute $h(B')$:
\begin{center}
\begin{tabular}{rrrrrrrrr}
 & $0$ & $-1$ & $-2$ & $-3$ & $0$ & $-1$ & $-2$ & $-3$ \\
$+$ & $0$ & $1$ & $2$ & $3$ & $4$ & $5$ & $6$ & $7$ \\
\hline
 & $0$ & $0$ & $0$ & $0$ & $4$ & $4$ & $4$ & $4$
\end{tabular}
\end{center}
Thus $B'$ is the board shown in Figure \ref{fig:fig5}.

\fig{.75}{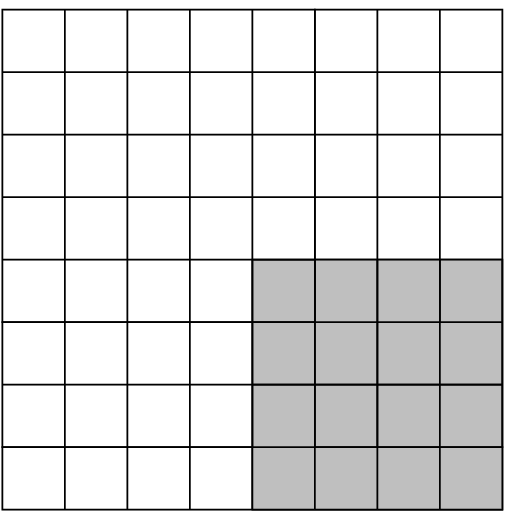}{The board $B^\prime$, rook-equivalent to $B$.}

In general, the board for even descents of permutations $\sigma \in S_{2n}$ has structure vector
\begin{equation*}
s(B) = \left( 0, 0, -1,-1, \ldots, -(n-1),-(n-1) \right),
\end{equation*}
and is thus rook-equivalent to the square $n \times n$ board $B^\prime$, which has structure vector
\begin{equation*}
s(B^\prime) = \left( 0, -1, \ldots, -(n-1), 0, -1, \ldots, -(n-1) \right).
\end{equation*}
The hit number $h_{s}$ for this square board $B^\prime$ inside the $2n \times 2n$ board is given by
\begin{equation*}
h_{s} = {n \choose s}^2 s! \cdot {n \choose n-s}(n-s)!\cdot n! = (n!)^2 {n \choose s}^2,
\end{equation*}
since we can first place $s$ rooks on the $n \times n$ board $B^\prime$ in ${n \choose s}^2 s!$ ways, then place $n-s$ rooks above $B^\prime$ in 
${n \choose n-s}(n-s)!$ ways, and finally place $n$ rooks in the left half of the $2n \times 2n$ board in $n!$ ways.

An important class of boards whose hit numbers have a simple product formula is the class of rectangular boards. 
One might ask whether there are other 
sets $X$ such that the board corresponding to $X$-descents
is rook equivalent to a rectangular board. In fact, Remark \ref{ap2} (in which $X = \{u+2, u+4, u+6,
\ldots, u+2m \}$) covers all possibilities. For example, consider the $2 \times 3$ rectangular board 
$B$ shown in Figure \ref{fig:fig18}.
\fig{.75}{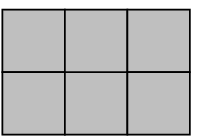}{The $2\times 3$ board $B$.}
We place this board in the lower right corner of an $8 \times 8$ board, as shown in Figure 
\ref{fig:fig6}.

\fig{.75}{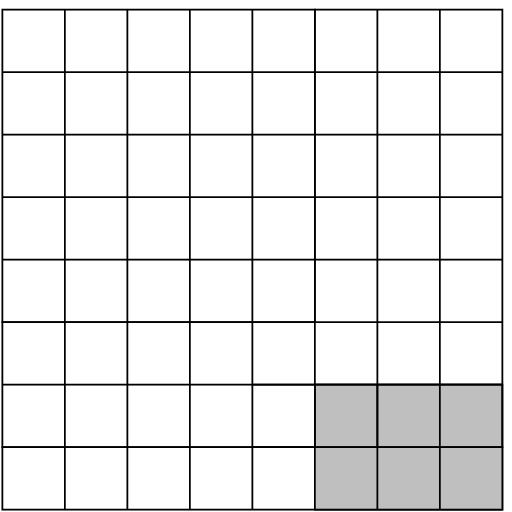}{The board $B$ placed in an $8 \times 8$ board.}

Note that for any set $X$ the board associated to $X$-descents has distinct rows, since
the row corresponding to $i \in X$ has length $i-1$. Rearranging the elements of $s(B)$ in weakly decreasing 
order gives the unique (see \cite{FoataSchutzenberger2}) board $B^\prime$ with distinct rows that is rook-equivalent to $B$.
In our example, we compute $s(B) = (0,-1,-2,-3,-4,-3,-4,-5)$. Thus $B$ is rook-equivalent to the
board $B'$ with structure vector $s(B') = (0,-1,-2,-3,-3,-4,-4,-5)$ and height vector $h(B') = (0,0,0,0,1,1,2,2)$. 
$B'$ is the board shown in Figure \ref{fig:fig7}.
\fig{.75}{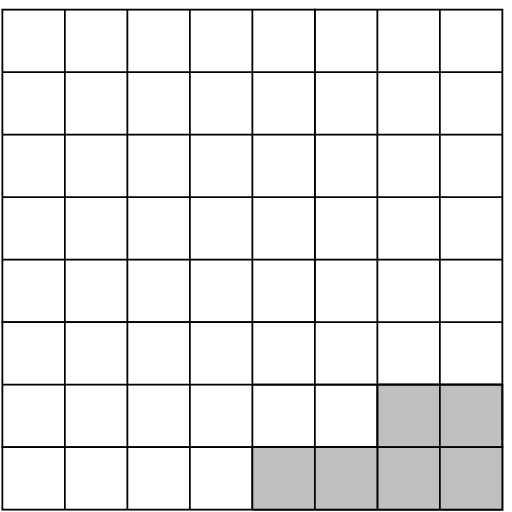}{The board $B^\prime$, rook-equivalent to $B$.}
Finally, to get the board for $X$-descents, 
we take the mirror image, and shift the rows upwards so that a row of length $i-1$ is in position $i$, 
as shown in Figure \ref{fig:fig8}.
\fig{.75}{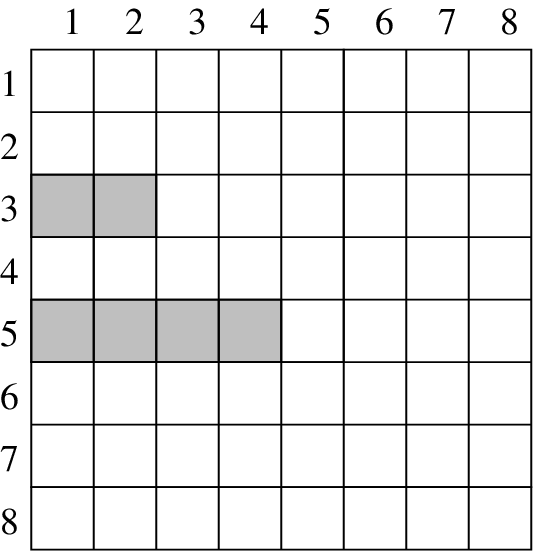}{The board $B^X_8$.}
Thus in our example, $X = \{3, 5\}$. In general, if we start with a rectangular $a \times b$ board ($a \leq b$), then the corresponding 
set is 
$X = \{u+2, u+4, u+6,
\ldots, u+2m \}$, where $m=a$ and $u=b-a$. For $n = 2m+u+v$, we get
\begin{eqnarray*}
P_{n,s}^X & = & {a \choose s}{b \choose s}s! \cdot {n-a \choose b-s}(b-s)! \cdot (n-b)! \\
& = & {m \choose s} {m+u \choose s}s! \cdot {m+u+v \choose m+u-s}(m+u-s)! \cdot (m+v)! \\
 & = & {m \choose s} {m+u+v \choose v+s}(m+u)!(m+v)!,
\end{eqnarray*}
exactly as in (\ref{eq:ap2}).

As noted in the introduction, we can also use the rook interpretation to translate the more general problem of counting $X,Y$-descents 
into one for which Corollaries 
\ref{comb1} and \ref{comb2} apply directly.
\begin{proposition}
Given subsets $X,Y \subseteq \mathbb{N}$ and a permutation $\sigma \in S_n$, let $B$ be the Ferrers board corresponding to the potential descent 
pairs $(i,j)$, where $i \in X_n$ and $j \in Y_n$. Let $B'$ be the unique Ferrers board rook-equivalent to $B$ that has distinct rows. Let $X' 
\subseteq [n]$ be the unique subset whose corresponding board (once empty rows and columns are deleted, and taking the 
mirror image) is $B'$. Then
\begin{equation*}
P_{n,s}^{X,Y} = P_{n,s}^{X'}.
\end{equation*}
\end{proposition}
\begin{proof}
By the previous discussion, we have $P_{n,s}^{X,Y} = h_s(B)$ and $P_{n,s}^{X^\prime} = h_s(B^\prime)$. But $B$ and $B^\prime$ are rook-equivalent,
and thus $h_s(B) = h_s(B^\prime)$ for all $s$.
\end{proof}

\begin{example}
Let $X = \{2,3,5,7,8\}, Y = \{1,2,4,5,6\}$, and $n=8$, so that the potential descent pairs are 
$21,31,32,51,52,54,71,72,74,75,76,81,82,84,85$, and $86$. Then 
$P_{8,s}^{X,Y} = h_s(B_8^{X,Y})$, where $B_8^{X,Y}$ is the board shown in Figure \ref{fig:fig9}.
\fig{.75}{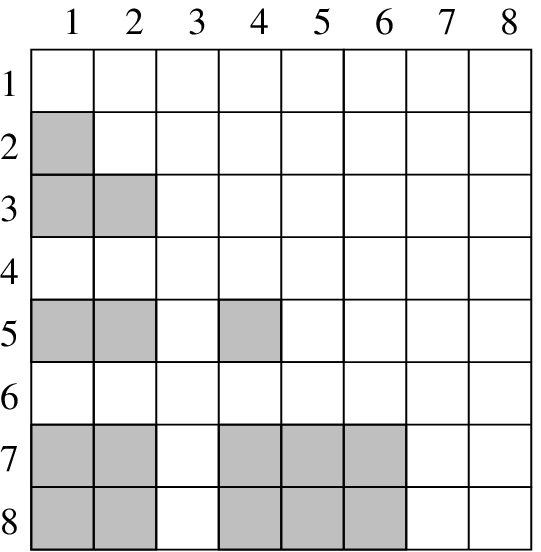}{The board $B_8^{X,Y}$.}
The unique board $B^\prime$ rook-equivalent to $B_8^{X,Y}$ that has distinct rows is shown in Figure 
\ref{fig:fig10}.
\fig{.75}{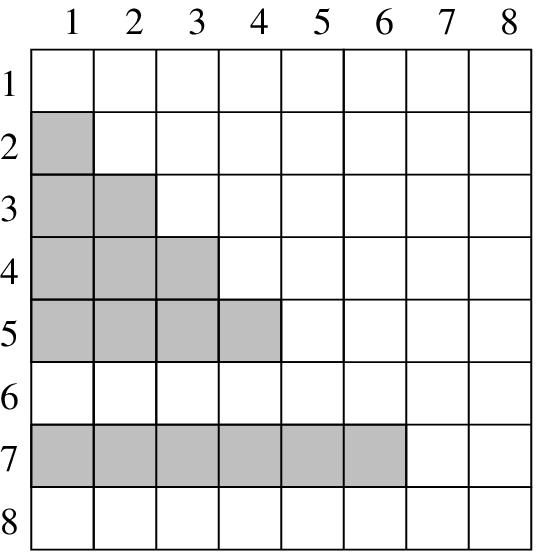}{The board $B^\prime$, rook-equivalent to $B_8^{X,Y}$.}
$B^\prime$ is the board for $X^\prime=\{2,3,4,5,7\}$, and so $P_{8,s}^{X,Y} = P_{8,s}^{X^\prime}$ for all $s$.
\end{example}


\section{Words}
\label{words-section}

Our results (\ref{eq:combXY1}) and (\ref{eq:combXY2}) extend easily to words. Let $\rho = (\rho_1,\rho_2,\ldots,\rho_m)$ be a 
composition of $n$, and let $R(\rho)$
be the rearrangement class of the word $1^{\rho_1}2^{\rho_2}\cdots m^{\rho_m}$ 
(i.e., $\rho_1$ copies of $1$, $\rho_2$ copies of $2$, etc.).
Given $X, Y \subseteq \mathbb{N}$, and a word $w \in R(\rho)$, define 
\begin{eqnarray*}
Des_{X,Y}(w) & = & \{ i : w_i > w_{i+1} 
~\&~ w_i \in X ~\&~ w_{i+1} \in Y \} ,\\
des_{X,Y}(w) & = & |Des_{X,Y}(w)|, \mbox{ and} \\
P_{\rho,s}^{X,Y} & = & \left| \{ w \in R(\rho) :des_{X,Y}(w) = s \} \right|.
\end{eqnarray*}
Then
\begin{theorem}\label{combword1}
\begin{equation}\label{eq:combword1}
P_{\rho,s}^{X,Y} = {a \choose \rho_{v_1},\rho_{v_2}, \ldots, \rho_{v_b}} \sum\limits_{r=0}^s (-1)^{s-r} {a+r \choose r}
{n+1 \choose s-r}
\prod\limits_{x\in X} {\rho_{x} + r + \alpha_{X,\rho,x}+\beta_{Y,\rho,x}  \choose \rho_{x} },
\end{equation}
where $X_m^c = \{v_1, v_2, \ldots, v_b \}, a = \sum\limits_{i=1}^b \rho_{v_i}$, and for any $x \in X_m$,
\begin{eqnarray*}
\alpha_{X,\rho,x} & = & \sum\limits_{\tiny \begin{array}{c} z \notin X \\ x < z \leq m \end{array} \normalsize}
\rho_z, \mbox{ and} \\ 
\beta_{Y,\rho,x} & = & \sum\limits_{\tiny \begin{array}{c} z \notin Y \\ 1 \leq z < x \end{array} \normalsize} \rho_z.
\end{eqnarray*}
\end{theorem}

\begin{proof}
We proceed as in the proof of Theorem \ref{combXY1}.
Fix $s$ and the composition $\rho = (\rho_1,\rho_2,\ldots,\rho_m)$ of $n$. Given $r$ such that $0 \leq r \leq s$, a 
$(\rho,s,r)^{X,Y}$-\emph{configuration} $c$ consists of an array of the elements of the 
multiset $\{ 1^{\rho_1}, 2^{\rho_2}, \ldots, m^{\rho_m} \}$, $r$ $+$'s, and $(s-r)$ $-$'s, satisfying
\begin{enumerate}
\item[(i)] each $-$ is either at the very beginning of the array or immediately follows a number, and
\item[(ii)] if $x$ and $y$ are consecutive numbers in the array, with $x \in X, y \in Y$, and $x>y$, i.e., if $(x,y)$ 
forms an $X,Y$-descent pair in the underlying word, then there must be at least one $+$ between $x$ and $y$.
\end{enumerate}
As an example, if $\rho = (2,3,1,4,2), X = \{2, 3, 5 \}$, and $Y=\{1, 3, 4\}$, the following is a $(\rho,5,3)^{X,Y}$-configuration. 
\begin{equation*}
c=4124-413-+25+42+5
\end{equation*}
In this example, the underlying word is $412441325425$. As before, we will let $c_1c_2 \cdots c_n$ denote the
underlying word of the $(\rho,s,r)^{X,Y}$-configuration $c$.

Let $C_{\rho,s,r}^{X,Y}$ be the set of all $(\rho,s,r)^{X,Y}$-configurations.  Then we claim 
that 
\begin{equation*}
\left| C_{\rho, s,r}^{X,Y} \right| = {a \choose \rho_{v_1},\rho_{v_2}, \ldots, \rho_{v_b}} {a+r \choose r}
{n+1 \choose s-r} \prod\limits_{x \in X} {\rho_{x} + r + \alpha_{X,\rho,x}+\beta_{Y,\rho,x} \choose \rho_{x} }.
\end{equation*}
That is, we can construct the set of $(\rho,s,r)^{X,Y}$-configurations as follows. First, we order the elements of the multiset $\{ v_1^{\rho_{v_1}}, \ldots, v_b^{\rho_{v_b}} \}$. This can be done in 
${a \choose \rho_{v_1},\rho_{v_2}, \ldots, \rho_{v_b}}$ ways. Next, we insert the  $r$ $+$'s. This can be done in ${a+r \choose r}$ ways.
Writing $X_m = \{ x_1 < x_2 < \cdots < x_{\left| X_m \right| } \}$, we can next place the elements of the 
multiset $\{ x_1^{\rho_{x_1}}, \ldots, x_{\left| X_m \right|}^{\rho_{x_{\left| X_m \right|}}} \}$ in 
$\prod\limits_{i=1}^{\left| X_m \right|} {\rho_{x_i} + r + \alpha_{X,\rho,x_i}+\beta_{Y,\rho,x_i} \choose \rho_{x_i}}$ 
ways, since after placing all copies of $x_1, x_2, \ldots, x_{i-1}$, the $\rho_{x_i}$ copies of $x_i$ can either go 
\begin{itemize}
\item[$\bullet$] immediately before any of the $\beta_{Y,\rho,x_i}$ elements of $\{1, 2, \ldots, x_{i-1} \}$ that is not in 
$Y$, or 
\item[$\bullet$] immediately before any of the $\alpha_{X,\rho,x_i}$ elements of $\{x_i+1,x_i+2, \ldots, m\}$ 
that is not in $X$, or
\item[$\bullet$] immediately before any of the $r$ $+$'s, or
\item[$\bullet$] at the very end of the array.
\end{itemize}
Thus, we have $1+r+\alpha_{X,\rho,x_i}+\beta_{Y,\rho,x_i} $ places in which we 
can insert letters equal to $x_i$. The number of ways in which we can place the $x_i$'s is therefore equal to the number of positive 
integral solutions of the equation $$z_1 + \cdots + z_{1+r + \alpha_{X,\rho,x_i}+\beta_{Y,\rho,x_i} } = \rho_{x_i},$$ which is well know to be 
${\rho_{x_i} + r + \alpha_{X,\rho,x_i}+\beta_{Y,\rho,x_i} \choose \rho_{x_i} }$.
Finally, we can place the $-$'s in ${n+1 \choose s-r}$ ways. 

As before, we define the weight $w(c)$ of an $(\rho,s,r)^{X,Y}$-configuration $c$ to be $-1$ to the number of $-$'s of $c$. It then follows that the RHS of 
(\ref{eq:combword1}) equals 
\begin{equation*}
\sum_{r=0}^s \sum_{c \in C_{\rho,s,r}^{X,Y}} w(c)
\end{equation*}
We now employ the identical involution $I$ as in the proof of Theorem \ref{combXY1}, this time on the set $C_{\rho, s}^{X,Y} 
= \bigsqcup\limits_{r=0}^s C_{\rho, s,r}^{X,Y}$. As before, we scan from left to right, and reverse the first 
sign that we can reverse without violating conditions (i) and (ii).
In our example above, the first place where we either encounter a sign that can be reversed is the $-$ after the second $4$.
Thus
\begin{equation*}
I(c)=4124+413-+25+42+5.
\end{equation*}
As was the case with the proof of Theorem \ref{combXY1}, it is simple to check that $I$ is a sign-reversing 
involution.

Now, suppose that $I(c)=c$. Then $c$ clearly can have no $-$'s, and so $r=s$ and $w(c) = 1$.
It must also be the case that no $+$'s can be reversed. Thus each of the $s$ $+$'s must occur singly in the 
middle of an $X,Y$-descent pair. It follows that the underlying word has exactly $s$ $X,Y$-descents.
\end{proof}

Taking $Y = \mathbb{N}$ in (\ref{eq:combword1}) gives the word analogue of Corollary \ref{comb1}.
\begin{corollary}\label{combword2}
\begin{equation}\label{eq:combword2}
P_{\rho,s}^X = {a \choose \rho_{v_1},\rho_{v_2}, \ldots, \rho_{v_b}} \sum\limits_{r=0}^s (-1)^{s-r} {a+r \choose r}
{n+1 \choose s-r}
\prod\limits_{x \in X} {\rho_{x} + r+ \alpha_{X,\rho,x} \choose \rho_{x}},
\end{equation}
where we write $P_{\rho,s}^X$ for $P_{\rho,s}^{X,\mathbb{N}}$.
\end{corollary}

\begin{corollary}\label{combword3}
Let $X = \{2 \}$ and $\rho = (a,b)$. Then
\begin{equation*}
P_{\rho,s}^X = {a \choose s}{b \choose s}. 
\end{equation*}
\end{corollary}
\begin{proof}
By (\ref{eq:combword2}) we have
\begin{eqnarray*}
P_{\rho,s}^X & = & {a \choose a} \sum\limits_{r=0}^s (-1)^{s-r} {a+r \choose r} {a+b+1 \choose s-r} {b+r \choose b} \\
 & = & \frac{(s+1)_a}{a!} \frac{(s+1)_b}{b!} \mbox{ }_3F_2 \left[ \begin{array}{rrr} -s, 
& -s, & -(a+b+1) \\ & -(a+s), & -(b+s) \end{array} \right] \\
 & = &  \frac{(s+1)_a}{a!} \frac{(s+1)_b}{b!} \frac{(a-s+1)_s (b-s+1)_s}{(a+1)_s (b+1)_s}\\
 & = & {a \choose s} {b \choose s}.
\end{eqnarray*}
\end{proof}

\begin{remark}
We can give a purely combinatorial proof of Corollary \ref{combword3} using Foata's transformation switching descents and excedences. The number
of rearrangements of $a$ $1$'s and $b$ $2$'s with exactly $s$ excedences is ${a \choose s}{b \choose s}$, since we have
to choose which $s$ of the first $a$ spots to be $2$'s (giving $s$ excedences) and which $s$ of the last $b$ spots to be $1$'s.
\end{remark}

\begin{corollary}\label{kwords}
Let $X = 2\mathbb{N}$ and $\rho = (\underbrace{k,k,\ldots,k}_{2n})$. Then
\begin{equation*}
P_{\rho,s}^X = {kn \choose k,k, \ldots, k}^2{kn \choose s}^2. 
\end{equation*}
\end{corollary}
\begin{proof}
By (\ref{eq:combword2}) we have
\begin{eqnarray*}
P_{\rho,s}^X & = & {kn \choose k,k, \ldots, k} \sum\limits_{r=0}^s (-1)^{s-r} {kn+r \choose r} {2kn+1 \choose s-r}
\prod_{i=1}^n {ki+r \choose k} \\
 & = & \frac{1}{(k!)^{2n}} \sum\limits_{r=0}^s (-1)^{s-r} (r+1)_{kn}^2 {2kn+1 \choose s-r} \\
 & = & \frac{(s+1)_{kn}^2}{(k!)^{2n}} \mbox{ }_3F_2 \left[ \begin{array}{rrr} -s, 
& -s, & -(2kn+1) \\ & -(kn+s), & -(kn+s) \end{array} \right] \\
 & = & {kn \choose k,k, \ldots, k}^2{kn \choose s}^2.
\end{eqnarray*}
\end{proof}

Finally, there is an alternative formula for $P_{\rho,s}^{X,Y}$. In addition to the notation of Theorem 
\ref{combword1}, let 
\begin{equation*}
\beta_{X,\rho,x} = \sum_{\tiny \begin{array}{c} z \notin X \\ 1 \leq z < x \end{array} \normalsize} \rho_z.
\end{equation*}
We have
\begin{theorem}\label{combword22}
\begin{equation}\label{eq:combword22}
P_{\rho,s}^{X,Y} = {a \choose \rho_{v_1},\rho_{v_2}, \ldots, \rho_{v_b}} \sum\limits_{r=0}^{n-a-s} (-1)^{n-a-s-r}
{a+r \choose r}{n+1 \choose n-a-s-r}
\prod_{x \in X} \binom{r+\beta_{X,\rho,x} -\beta_{Y,\rho,x}}{\rho_{x}},
\end{equation}
where we use the convention that $\binom{p}{q} = 0$ if $p < 0$.
\end{theorem}

\begin{proof}
The proof is analogous to that of Theorem \ref{combXY2}. Given a fixed composition $\rho = (\rho_1,\rho_2,\ldots,\rho_m)$ of $n$, and $s 
\geq r \geq 0$, a $\overline{(\rho,s,r)}^{X,Y}$-\emph{configuration} is an array of
the elements of the multiset $\{ 1^{\rho_1}, 2^{\rho_2}, \ldots, m^{\rho_m} \}$, together with $r$ $+$'s and
$(n-a-s-r)$ $-$'s, satisfying
\begin{enumerate}
\item[(i)] each $-$ is either at the very beginning of the string or immediately follows a number,
\item[(ii)] if $c_i \in X, 1 \leq i < n$, and $(c_i, c_{i+1})$ is not an $X,Y$-descent pair of the underlying 
word, 
then there must be at least one $+$ between $c_i$ and $c_{i+1}$, and
\item[(iii)] if $c_n \in X$, then $c_n$ must be followed by at least one $+$.
\end{enumerate}
As an example, if $X = \{2, 3, 6 \},Y = \{1,2,5\}$, and $\rho = 
(2,1,3,2,1,1)$, then the following is a 
$\overline{(\rho,1,2)}^{X,Y}$-configuration. 
\begin{equation*}
213+43+3-165-4
\end{equation*}

Let $\overline{C}_{\rho,s,r}^{X,Y}$ be the set of all $\overline{(\rho,s,r)}^{X,Y}$-configurations.  Then we claim 
that 
\begin{equation*}
\left| \overline{C}_{\rho, s,r}^{X,Y} \right| = {a+r \choose r} \binom{a}{\rho_{v_1},\rho_{v_2}, \ldots, \rho_{v_b}}  {n+1 \choose n-a 
-s-r} \prod_{x \in X} 
\binom{r+\beta_{X,\rho,x} -\beta_{Y,\rho,x}}{\rho_{x}}.
\end{equation*} 
That is, we can construct the set of $\overline{(\rho,s,r)}^{X,Y}$-configurations as follows. First, we order the elements of the multiset $\{ v_1^{\rho_{v_1}}, \ldots, v_b^{\rho_{v_b}} \}$. This can be done in 
${a \choose \rho_{v_1},\rho_{v_2}, \ldots, \rho_{v_b}}$ ways. Next, we insert the  $r$ $+$'s. This can be done in ${a+r \choose r}$ ways.
Next, consider the choices for placing the elements of the multiset $\{ x_1^{\rho_{x_1}}, x_2^{\rho_{x_2}}, \ldots, 
x_{ \left| X_m^c \right| }^{\rho_{x_{\left| X_m^c \right|}}}\}$ in the order $x_1 < x_2 < \cdots < x_{\left| X_m^c \right|}$. First,
there are $r+\beta_{X,\rho,x_1}-\beta_{Y,\rho,x_1}$ spaces in which to insert the $x_1$'s, since $x_1$ can either go 
immediately before of any $+$, or immediately before any element of $Y$ which is 
less than $x_1$. Note that unlike the situation in Theorem \ref{combword1}, no more than one copy of $x_1$ can go in 
any particular available space. Thus there are $\binom{r+\beta_{X,\rho,x_1}-\beta_{Y,\rho,x_1}}{\rho_{x_1}}$ ways to 
place the $x_1$'s. 

In general, having placed all copies of 
$x_1, x_2, \ldots, x_{i-1}$, we cannot place $x_i$ immediately before 
some $y \in Y, y < x_i$, which earlier had an 
element of the multiset $\{ x_1^{\rho_{x_1}}, x_2^{\rho_{x_2}}, \ldots,  x_{i-1}^{\rho_{x_{i-1}}}\}$ placed 
immediately before it. 
Similarly, we cannot place $x_i$ immediately before any $+$ which earlier had an element of the multiset 
$\{ x_1^{\rho_{x_1}}, x_2^{\rho_{x_2}}, \ldots,  x_{i-1}^{\rho_{x_{i-1}}}\}$ placed immediately before it. It then 
follows that
there are
\begin{equation*}
r+\sum_{1\leq z<x_i} \rho_z - \beta_{Y,\rho,x_i} - \sum_{j=1}^{i-1} \rho_{x_j} = r+\beta_{X,\rho,x_i} -\beta_{Y,\rho,x_i}
\end{equation*}
spaces in which to insert the $x_i$'s. Thus, there are total of $\prod_{i=1}^{\left| X_m^c \right|} 
\binom{r+ \beta_{X,\rho,x_i} -\beta_{Y,\rho,x_i}}{\rho_{x_i}}$ ways to place 
all of the copies of $x_1,x_2, \ldots,x_{\left| X_m^c \right|}$, given our placement of the copies of $v_1, v_2, \ldots, v_b$.
Finally, we can place the $-$'s in ${n+1 \choose n-a-s-r}$ ways. 

The remainder of the proof follows exactly as in Theorem \ref{combXY2}. The weight of a configuration is defined in the same way, and 
applying the same sign-reversing involution $I$ cancels out all configurations except those corresponding to words with exactly $s$ 
$X,Y$-descent pairs.
\end{proof}

In the special case $Y= \mathbb{N}$, Theorem \ref{combword22} reduces to the following. 
\begin{corollary}\label{combword221}
\begin{equation*}\label{eq:combword221}
P_{\rho,s}^X = {a \choose \rho_{v_1},\rho_{v_2}, \ldots, \rho_{v_b}} \sum\limits_{r=0}^{n-a-s} (-1)^{n-a-s-r} {a+r \choose r}
{n+1 \choose s-r}
\prod\limits_{x \in X} \binom{r+\beta_{X,\rho,x}}{\rho_{x}}.
\end{equation*}
\end{corollary}



\section{$X,Y,Z$-descents}

As mentioned in the Introduction, a more general problem is to study the class of polynomials 
\begin{equation*}
P_n^{X,Y,Z}(x) = \sum_{s \geq 0} P_{n,s}^{X,Y,Z}x^s := \sum_{\sg \in S_n} x^{des_{X,Y,Z}(\sg)},
\end{equation*}
where for any subsets $X$, $Y$, and $Z$ of $\mathbb{N}$, and 
permutation $\sg = \sg_1 \sg_2 \cdots \sg_n \in S_n$,
\begin{eqnarray*}\label{DESXYZ}
Des_{X,Y,Z}(\sg) & = & \{i:\sg_i > \sg_{i+1} \ \& \ \sg_i \in X, \sg_{i+1} \in 
Y \ \& \ \sg_i - \sg_{i+1} \in Z\}, \mbox{ and} \\
des_{X,Y,Z}(\sg) & = & |Des_{X,Y,Z}(\sg)|.
\end{eqnarray*}
We say shall that $(\sg_i,\sg_{i+1})$ is an $X,Y,Z$-\emph{descent} if 
$i \in Des_{X,Y,Z}(\sg)$. 
The polynomials studied in this paper are thus the special case $Z = \mathbb{N}$ of the polynomials 
$P_n^{X,Y,Z}(x)$.

In many cases, we can obtain formulas for the coefficients $P_{n,s}^{X,Y,Z}$ from our previous formulas. That is, in many cases, 
the possible  $X,Y,Z$-descents under Foata's transformation corresponds to a board which is rook equivalent to a Ferrers board. In such cases, we can 
use the formulas for hit polynomials or our formulas for $P_{n,s}^{X,Y}$ or 
$P_{n,s}^{X}$ to obtain formulas for $P_{n,s}^{X,Y,Z}$. For example, 
let $E = 2\mathbb{N}$ and $\mathbb{N}_{\geq k} = \{k,k+1,k+2, \ldots \}$.  
Under Foata's 
transformation, the ${E,E,\mathbb{N}_{\geq 2k}}$-descents correspond to 
excedences of the form $\begin{array}{c} 2s \\ 2t \end{array}$ where $2t -2s \geq 2k$. For example, 
if $k =2$ and $n =20$, then we would consider the board $B_{20}^{E,E,\mathbb{N}_{\geq 4}}$ shown in Figure 
\ref{fig:fig11}. 
\fig{.5}{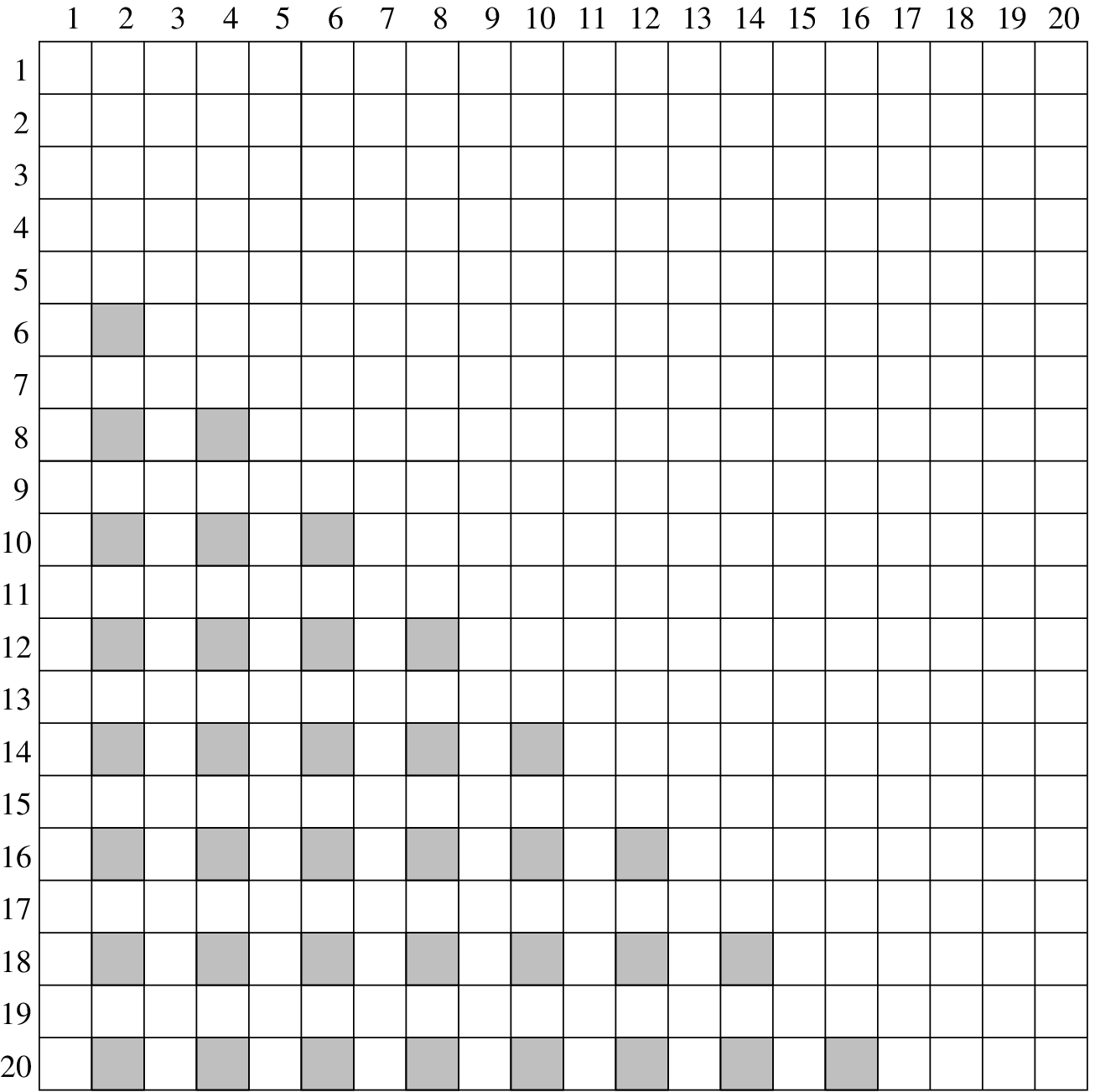}{The board $B_{20}^{E,E,\mathbb{N}_{\geq 4}}$.}
It is then easy to see that $B_{20}^{E,E,\mathbb{N}_{\geq 4}}$ is rook equivalent to the board 
$B_{20}^{\{2,3,4,5,6,7,8,9\}}$ pictured in Figure \ref{fig:fig12}, which is the board we would consider when computing 
the polynomial $P_{20}^{\{2,3,4,5,6,7,8,9\}}(x)$. It follows that 
$P_{20}^{E,E,\mathbb{N}_{\geq 4}}(x)= P_{20}^{\{2,3,4,5,6,7,8,9\}}(x)$.
Thus we can use Corollary \ref{comb1} or Corollary \ref{comb2} 
to give explicit 
formulas for $P_{20,s}^{E,E,\mathbb{N}_{\geq 4}} = 
P_{20,s}^{\{2,3,4,5,6,7,8,9\}}$.

\fig{.5}{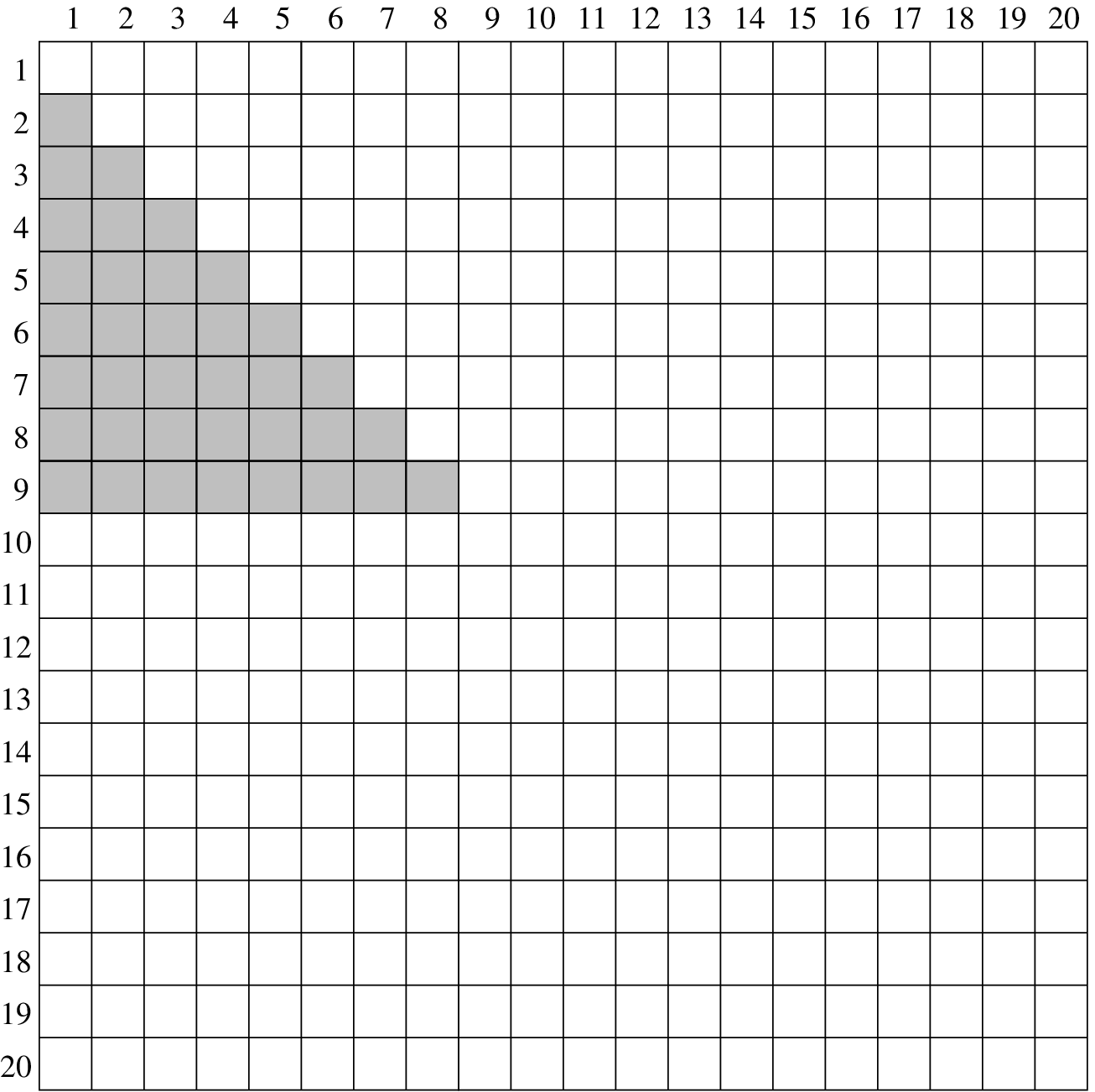}{The board $B_{20}^{\{2,3,4,5,6,7,8,9\}}$.}

The problem of computing $P_{n,s}^{X,Y,Z}$ for arbitrary sets $X, Y$, and $Z$ seems to be difficult in large part because the 
board corresponding to $X,Y,Z$-descents is not a Ferrers board in general.
This can be seen in Example \ref{XYZrook}, in which $X=E,Y=O$, and $Z = \{1, 3\}$. In some special cases, however, the rook-placement 
formulation will enable us to derive a formula
for $P_{n,s}^{X,Y,Z}$ involving a double or triple sum. One such example is given below.

\begin{example}
Let $X  = \{ z : z = 4, 5, \mbox{ or } 6 \mod 6 \}, Y = \{ z : z = 1, 2, \mbox{ or } 3 \mod 6 \}$, and $Z = \{ 1, 2, 3, 4, 5, 6 \}$. For 
$n=12$, the board $B_{12}^{X,Y,Z}$ corresponding to $X,Y,Z$-descents is shown in Figure \ref{fig:fig13}.
\fig{.5}{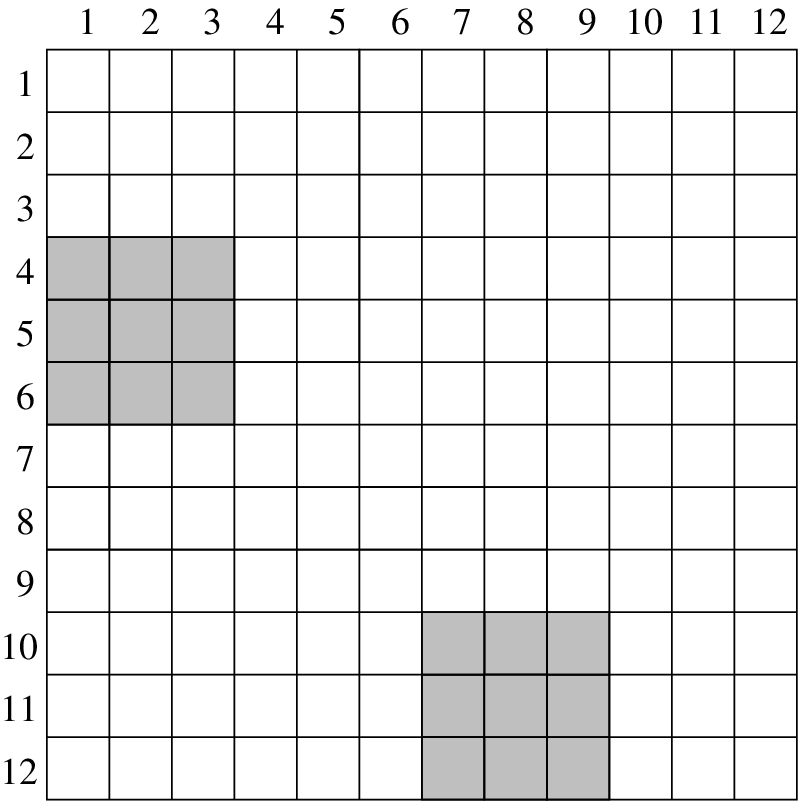}{The board $B_{12}^{X,Y,Z}$.}
By permuting rows and columns we see that $B_{12}^{X,Y,Z}$ is rook equivalent to the board $B^\prime$ shown in 
Figure \ref{fig:fig14}.
\fig{.5}{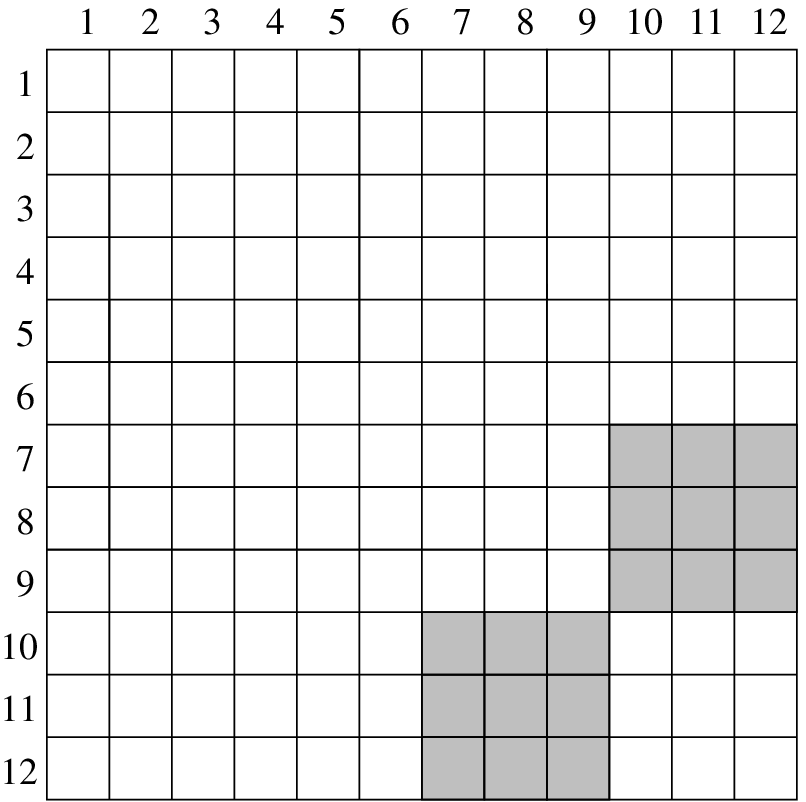}{The board $B^\prime$, rook-equivalent to $B_{12}^{X,Y,Z}$.}
To compute the hit number $h_s$ of $B^\prime$, we think of placing the rooks in several successive steps, as 
indicated by the numbers on the diagram in Figure \ref{fig:fig15}. 

\fig{.5}{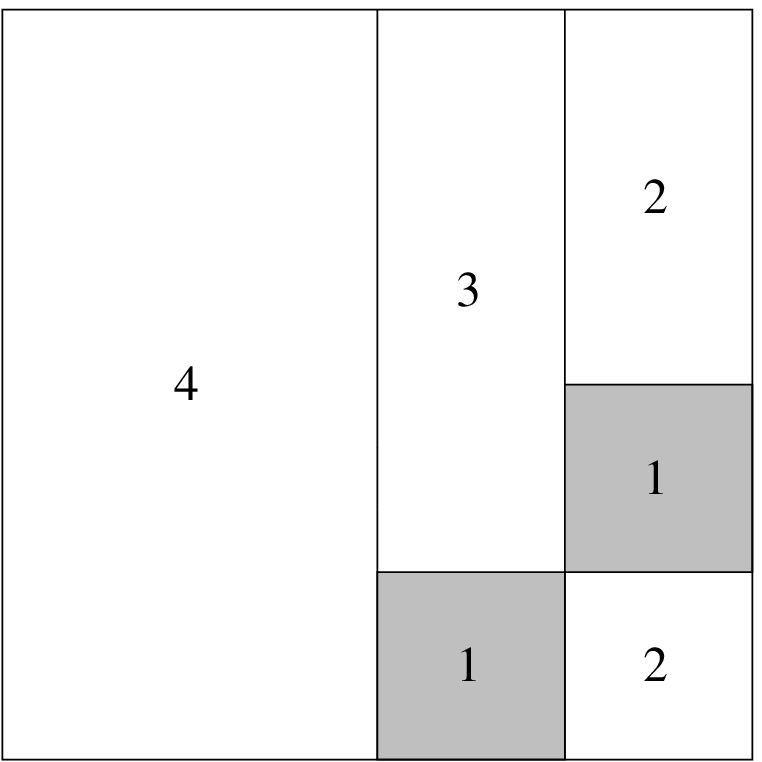}{A diagram indicating the order of placement of rooks on $B^\prime$.}

Suppose we first place $p$ rooks on the lower-left $3 \times 3$ 
block of $B^\prime$ and $q$ rooks on the upper-right $3 \times 3$ block of $B^\prime$, where $p + q =s$. This first step can be done in 
${3 \choose p} {3 \choose p}p!
\cdot {3 \choose q}{3 \choose q} q!$ ways. Next, we place a total of $3-q$ rooks in the regions marked `2'; this 
can be done in ${9-p \choose 3-q} 
(3-q)!$ ways. We then place a total of $3-p$ rooks in the regions marked `3'; this can be done in ${9-p \choose 3-p}(3-p)!$ ways. At this point, there
are six rooks left to place, one in each of the six columns of the regions marked `4'. Since six rooks have already been placed, only six rows remain open.
Hence there are $6!$ ways to do this last step. Thus we have
\begin{eqnarray*}
P_{12,s}^{X,Y,Z} & = & \sum\limits_{\begin{array}{c} p+q = s \\ 0 \leq p \leq 3 \\ 0 \leq q \leq 3 \end{array}} {3 \choose p} {3 \choose p}p! 
\cdot {3 \choose q}{3 \choose q} q! \cdot {9-p \choose 3-q} (3-q)! \cdot {9-p \choose 3-p}(3-p)! \cdot 6!\\
 & = & \sum\limits_{\begin{array}{c} p+q = s \\ 0 \leq p \leq 3 \\ 0 \leq q \leq 3 \end{array}} \frac{(3!)^4}{((3-p)!)^2((3-q)!)^2 p! q!} \frac{((9-p)!)^2}{(6-p+q)!}.
\end{eqnarray*}
\end{example}
\ \\

\section{Further Questions}
\label{questions-section}

In this section, we discuss some open questions and directions for further 
research.

\noindent {\bf Relations between $P_{\rho,s}^X$ and $P_{n,s}^X$.} It is easy to see that for
any composition $\rho = (\rho_1, \rho_2, \ldots, \rho_m)$ of $n$, we have
\begin{equation*}
\rho_1! \rho_2! \cdots \rho_m! |R(\rho)| = |S_n|.
\end{equation*}
In fact there is a natural bijection
\begin{eqnarray*}
\chi: (S_{\rho_1} \times S_{\rho_2} \times \cdots \times S_{\rho_m}) \times R(\rho) & \longrightarrow & S_n, \\
 ((\phi^{(1)},\phi^{(2)}, \ldots, \phi^{(m)}),w) & \longmapsto & \sigma,
\end{eqnarray*}
where $\sigma$ is obtained from $w$ by replacing the $i$th occurrence of $1$ by $\phi^{(1)}_i$, the $i$th occurrence of $2$ by $\rho_1 + \phi^{(2)}_i$, and so on.
For example, $\chi((21,312),12212) = 25314$.
As another example, if $\phi^{(j)}$ is the identity permutation for each $j$, then $\chi((\phi^{(1)},\phi^{(2)}, \ldots, \phi^{(m)}),w)$ is just the
usual \emph{standardization} of $w$, written $std(w)$, which is obtained from $w$ by replacing the $i$th occurrence of $1$ with $i$, 
the $i$th occurrence of $2$ with $\rho_1 + i$, and so on. However, in general it is \emph{not} true that
\begin{equation}\label{eq:std}
\rho_1! \rho_2! \cdots \rho_m! P_{\rho,s}^X = P_{n,s}^X,
\end{equation}
for arbitrary $\rho$ and $X$. For example, if $X = \mathbb{N}$ and we are counting descents without restriction, then $P_{\rho,s}^X = 0$ 
for
$s > \sum_{i=1}^{m-1} \rho_i$, since the largest number $m$ cannot be the bottom of a descent. 
On the other hand, $P_{n,s}^X$ is non-zero for all $0 \leq s \leq n-1$.

In light of the above, the following consequence of
Corollaries \ref{ap1} and \ref{kwords} is quite surprising. For $X = 2\mathbb{N}$ and $\rho = 
(\underbrace{k,k,\ldots,k}_{2n})$, we have
\begin{equation*}
 (k!)^{2n}P_{\rho,s}^X = P_{2kn,s}^X.
\end{equation*}
This identity does not follow by applying the bijection $\chi$. For example, $w = 41421323$ has two even descents, while $std(w) = 71832546$ has one. 
We therefore ask for an explicit bijection
\begin{equation*}
\psi: (\underbrace{S_k \times S_k \times \cdots \times S_k}_{2n}) \times R(\underbrace{k,k,\ldots,k}_{2n})  \longrightarrow  S_{2kn},
\end{equation*}
satisfying
\begin{equation*}
\overleftarrow{des}_E(\psi(w)) = \overleftarrow{des}_E(w)
\end{equation*}
for all $w \in R(\underbrace{k,k,\ldots,k}_{2n})$. We also ask if there are other sets $X$ and compositions $\rho$ for which 
(\ref{eq:std}) holds.
\ \\

\noindent {\bf $q$-analogues.} The following $q$-analogue of the numbers $P^{X}_{n,s} = P^{X,\mathbb{N}}_{n,s}$ exists.
Let 
\begin{equation*}
[n]_q = 1+q+q^2+\ldots+q^{n-1}.
\end{equation*}
Let 
$\Delta_{n+1}^q$ and $\Gamma_{n+1}^q$ be the operators defined as
\begin{eqnarray*}
\Delta_{n+1} & : & x^s \longrightarrow [s]_q x^{s-1} + q^s[n+1-s]_q x^s 
\\
\Gamma_{n+1} & : & x^s \longrightarrow [s+1]_q x^s + q^{s+1}[n-s]_q x^{s+1}. 
\end{eqnarray*}
Given a subset $X \subseteq\mathbb{N}$, we define the
polynomials $P_n^X(q,x)$ by $P_0^X(q,x)=1$, and
\begin{equation*}
P_{n+1}^X(q,x) = \left\{ \begin{array}{ll} \Delta_{n+1}^q (P_n^X(q,x)), & \mbox{ if } n+1 \not\in X, \mbox{ and} 
\\
\Gamma_{n+1}^q (P_n^X(q,x)), & \mbox{ if } n+1 \in X. \end{array} \right.
\end{equation*}
We define the coefficient polynomials $P_{n,s}^X(q)$ by setting $P_n^X(q,x) = \sum\limits_{s \geq 
0}P_{n,s}^X(q)x^s$. The polynomials
$P_{n,s}^X(q)$ satisfy a recursion analogous to (\ref{XYrec}), and $q$-analogues of the formulas (\ref{eq:comb1}) and (\ref{eq:comb2})
have been found. These results, along with a combinatorial interpretation of a Mahonian statistic $stat$ satisfying
\begin{equation*}
P_{n,s}^X(q) = \sum\limits_{\tiny \begin{array}{c} \sigma \in S_n, \\ des_{X}(\sigma) = s \end{array}} q^{stat(\sigma) \normalsize},
\end{equation*}
will be presented in an upcoming paper by the authors and J. Liese.

\ \\

\noindent {\bf Pattern matchings.} One can put our results in a more general context of pattern matchings 
in permutations as follows. Given any sequence $\sg = \sg_1 \sg_2 \cdots \sg_n$ of distinct integers,
we let $red(\sg)$ be the permutation that results by replacing the
$i$-th smallest integer that appears in the sequence $\sg$ by $i$.
For example, if $\sg = 2~7~5~4$, then $red(\sg) = 1~4~3~2$. Given a
permutation $\tau$ in the symmetric group $S_j$, we define a
permutation $\sg = \sg_1 \sg_2 \cdots \sg_n \in S_n$ to have a
$\tau$-match at place $i$ provided $red(\sg_i \sg_{i+1} \cdots \sg_{i+j-1}) =
\tau$.  Let $\tau\mbox{-}mch(\sigma)$ be the number of
$\tau$-matches in the permutation $\sg$. To prevent confusion, we
note that a permutation not having a $\tau$-match is different than
a permutation being $\tau$-avoiding.  A permutation is called
$\tau$-\emph{avoiding} if there are no indices $i_1 < i_2 < \cdots < i_j$ such
that $ red(\sg_{i_1} \sg_{i_2} \cdots \sg_{i_j}) = \tau$. For example, if
$\tau = 2~1~4~3$, then the permutation $3~2~1~4~6~5$ does not have a
$\tau$-match but it does not avoid $\tau$ since $red(2~1~6~5) =
\tau$. In the case where $|\tau|=2$, $\tau\mbox{-}mch(\sigma)$ reduces
to familiar permutation statistics. That is, if $\sg = \sg_1 \sg_2 \cdots
\sg_n \in S_n$, let $Des(\sg) = \{i:\sg_i > \sg_{i+1}\}$ and
$Rise(\sg) = \{i:\sg_i < \sg_{i+1}\}$. Then it is easy to see that
$(2~1)\mbox{-}mch(\sg) = des(\sg) = |Des(\sg)|$ and
$(1~2)\mbox{-}mch(\sg) = rise(\sg) = |Rise(\sg)|$.
A number of recent publications have analyzed the distribution of
$\tau$-matches in permutations. See, for example,
\cite{Eliz,Kit1,Kit2}.

We can consider a more refined pattern-matching condition
where we take into account conditions involving equivalence mod $k$
for some integer $k \geq 2$. That is, suppose we fix $k \geq 2$ and
we are given some sequence of distinct integers $\tau = \tau_1 \tau_2
\cdots \tau_j$. Then we say that a permutation $\sg = \sg_1 \sg_2 \cdots
\sg_n \in S_n$ has a $\tau$-$k$-equivalence match at place $i$
provided $red(\sg_i \sg_{i+1} \cdots \sg_{i+j-1}) = red(\tau)$ and for all $s
\in \{0, 1, \ldots, j-1\}$, $\sg_{i+s} = \tau_{1+s} \mod k$. For
example, if $\tau = 1~2$ and $\sg = 5~1~7~4~3~6~8~2$, then $\sg$ has
$\tau$-matches starting at positions 2, 5, and 6. However, if $k=2$,
then only the $\tau$-match starting at position 5 is a
$\tau$-$2$-equivalence match. (Later, it will be explained that the
$\tau$-match starting at position 2 is a $(1~3)$-$2$-equivalence
match and the $\tau$-match starting a position 6 is a
$(2~4)$-$2$-equivalence match.) Let
$\tau\mbox{-}k\mbox{-}emch(\sigma)$ be the number of
$\tau$-$k$-equivalence matches in the permutation $\sg$.   

More generally, if $\Upsilon$ is a set of sequences of distinct
integers of length $j$, then we say that a permutation $\sg = \sg_1
\sg_2 \cdots \sg_n \in S_n$ has a $\Upsilon$-$k$-equivalence match at
place $i$ provided that there is a $\tau \in \Upsilon$ such that
$red(\sg_i  \sg_{i+1} \cdots \sg_{i+j-1}) = red(\tau)$ and for all $s \in \{0,
\ldots, j-1\}$, $\sg_{i+s} = \tau_{1+s} \mod k$. Let
$\Upsilon\mbox{-}k\mbox{-}emch(\sigma)$ be the number of
$\Upsilon$-$k$-equivalence matches in the permutation $\sg$.

One can then study the polynomials
\begin{eqnarray*}
T_{\tau,k,n}(x) &=& \sum_{\sg \in S_n}
x^{\tau\mbox{-}k\mbox{-}emch(\sigma)} =
\sum_{s=0}^n T_{\tau,k,n}^s x^s \ \mbox{and} \label{eq:k-tau} \\
U_{\Upsilon,k,n}(x) &=& \sum_{\sg \in S_n}
x^{\Upsilon\mbox{-}k\mbox{-}emch(\sigma)} =
\sum_{s=0}^n U_{\Upsilon,k,n}^s x^s.\label{eq:k-Up}
\end{eqnarray*}
In particular, suppose that we focus on the special cases of these
polynomials where we consider only patterns of length 2. That is,
fix $k \geq 2$ and let $A_k$ equal the set of all sequences $(a~b)$
such that $1 \leq a < b \leq 2k$ and there is no lexicographically
smaller sequence $x~y$ having the property that $x \equiv a \mod k$
and $y \equiv b \mod k$. For example, $$A_4 = \{1~2, 1~3, 1~4, 1~5,
2~3, 2~4, 2~5, 2~6, 3~4, 3~5, 3~6, 3~7, 4~5, 4~6, 4~7, 4~8\}.$$ Let
$D_k = \{b~a: a~b \in A_k\}$ and $E_k = A_k \cup D_k$. Thus $E_k$
consists of all $k$-equivalence patterns of length 2 that we could
possibly consider. Note that if $\Upsilon = A_k$, then
$\Upsilon\mbox{-}k\mbox{-}emch(\sigma) = rise(\sg)$ and if $\Upsilon
= D_k$, then $\Upsilon\mbox{-}k\mbox{-}emch(\sigma) = des(\sg)$.

Liese \cite{Liese} studied the polynomials $U_{\Upsilon,k,n}^s$ where 
$\Upsilon$ consists of patterns of length 2. For example, he showed that one 
can use inclusion-exclusion to find a formula for
$U_{\Upsilon,k,n}^s$ for any $\Upsilon \subset E_k$ in terms of
certain rook numbers of a sequences of boards associated with
$\Upsilon$. The same is true for coefficients of the 
the polynomials $P_n^{X,Y,Z}(x)$. This approach leads to completely 
different formulas than the ones produced in this paper. 
While this approach is straightforward, it is
unsatisfactory since it reduces the computation of
$U_{\Upsilon,k,n}^s$ to another difficult problem, namely, computing
rook numbers for general boards. 

Liese \cite{Liese} was able to give direct formulas for the coefficients $T_{\tau,k,n}^s$ where $\tau
\in E_k$. For example, in the case where $\tau= (1~k)$, his results
imply that for all $0 \leq s \leq n$ and for all $0 \leq j \leq
k-1$,
\begin{eqnarray*}
&& T_{(1~k),k,kn+j}^s = \\
&& \hspace{0.2cm} ((k-1)n+j)!\sum_{r=0}^s (-1)^{s-r} ((k-1)n+j+r)^n
\binom{(k-1)n+j+r}{r} \binom{kn+j+1}{s-r}, \mbox{ and} \\
&& T_{(1~k),k,kn+j}^s = \\
&& \hspace{0.2cm} ((k-1)n+j)!\sum_{r=0}^{n-s} (-1)^{n-s-r} (1+r)^n
\binom{(k-1)n+j+r}{r} \binom{kn+j+1}{n-s-r}.
\end{eqnarray*}
These two formulas are easily seen to be equivalent to special cases of our formulas. However, Liese 
has produced explicit formulas
for $U_{\Upsilon,k,n}^s$ in the special case where $\Upsilon$ is a
subset of the form $\{(x_1,y_1),(x_2,y_2),\ldots,(x_n,y_n)\}$, where
for all $i,j$ $y_i \equiv y_j \mod k$ and either $\Upsilon \subseteq
A_k$ or $\Upsilon \subseteq D_k$. These formulas cannot always be reduced to special cases of our formulas.



\begin{thebibliography}{99}

\bibitem{Comtet}
L.\ Comtet,
Permutations by Number of Rises; Eulerian Numbers,
in ``Advanced Combinatorics: The Art of Finite and Infinite Expansions,''
rev. enl. ed. Dordrecht, Netherlands: Reidel, 1974. 


\bibitem{Eliz}
S.\ Elizalde, M.\ Noy, Consecutive patterns in
permutations, \emph{Adv. in Appl. Math.} \textbf{30} (2003), no.~1-2, 110--125, Formal power series and
  algebraic combinatorics (Scottsdale, AZ, 2001).

\bibitem{EhrenborgSteingrimsson}
R.\ Ehrenborg, E.\ Steingrimsson,
The Excedence Set of a Permutation,
\emph{Adv. in Appl. Math.},
{\bf 24} (2000), 284--299.

\bibitem{FoataSchutzenberger1}
D.\ C.\ Foata, M.\ P.\ Sch\"{u}tzengerger,
``Theorie Geometriques des Polynomes Euleriens,''
Lecture Notes in Math., {\bf 138}, Springer-Verlag, Berlin, 1970.

\bibitem{FoataSchutzenberger2}
D.\ C.\ Foata, M.\ P.\ Sch\"{u}tzengerger,
On the rook polynomials of Ferrers relations,
in ``Combinatorial Theory and Its Applications, 2,''
Colloq. Math. Janos Bolyai, Vol. 4, pp. 413--436, North-Holland, Amsterdam, 1970.

\bibitem{Gasper}
G.\ Gasper,
Summation Formulas for Basic Hypergeometric Series,
\emph{Siam. J. Math. Anal.},
{\bf 12} (1981), 196--200.

\bibitem{GasperRahman}
G.\ Gasper, M.\ Rahman,
Basic Hypergeometric Series, in ``Encyclopedia of Math. and its Applications,''
Cambridge Univ. Press,
Cambridge, MA, 1990.

\bibitem{GoldmanJoichiWhite}
J.\ R.\ Goldman, J.\ T.\ Joichi, D.\ E.\ White,
Rook theory. I. Rook equivalence of Ferrers boards,
\emph{Proc. Amer. Math. Soc.},
{\bf 52} (1975), 485--492.

\bibitem{Haglund}
J.\ Haglund,
Rook Theory and Hypergeometric Series,
\emph{Adv. in App. Math.},
{\bf 17} (1996), 408--459.

\bibitem{HallLieseRemmel}
J.\ Hall, J.\ Liese, J.\ Remmel,
$q$-analogues of formulas counting descent pairs with prescribed tops and bottoms,
in preparation.

\bibitem{KaplanskyRiordan}
I.\ Kaplansky, J.\ Riordan,
The problem of the rooks and its applications,
\emph{Duke Math. J.},
{\bf 13} (1946), 259--268.

\bibitem{Kit2}
S.\ Kitaev, Partially ordered generalized patterns, to
appear in Discrete Math.

\bibitem{Kit1}
S.\ Kitaev, Generalized patterns in words and
permutations, Ph.D.
  thesis, Chalmers University of Technology and G\"{o}teborg University, 2003.

\bibitem{Kit3}
S.\ Kitaev and T.\ Mansour, Partially ordered
generalized patterns and $k$-ary words, \emph{Ann. Comb.} \textbf{7} (2003), no.~2, 191--200.

\bibitem{KitaevRemmel1}
S.\ Kitaev, J.\ Remmel,
Classifying Descents According to Parity,
{\tt math.CO/0508570} 

\bibitem{KitaevRemmel2}
S.\ Kitaev, J.\ Remmel,
Classifying Descents According to Equivalence mod $k$,
{\tt math.CO/0604455}

\bibitem{Liese} J. Liese, Classifying ascents and descents with specified equivalences mod $k$, Proceedings of the 
2006 International Conference on Formal Power Series and Algebraic Combinatorics.

\bibitem{LieseRemmel}
J.\ Liese, J.\ Remmel,
$q$-Analogues of formulas for the number of ascents and descents with specified equivalences mod $k$,
\emph{Permutation Patterns}, 2006

\end{thebibliography}
\end{document}